\def\DateTime{5/January/2007, 17:30(JP)}
\def\Version{Version $2.0$}
\def\yes{\if00}
\def\no{\if01}
\def\iftwelvept{\no}

\def\ifusepdf{\no}
\def\ifpsfont{\yes}

\iftwelvept
\documentclass[leqno,12pt]{amsart}
\else
\documentclass[leqno]{amsart}
\fi
\usepackage{amssymb}
\usepackage{amscd}
\usepackage{latexsym}
\usepackage{verbatim}
\usepackage[all]{xy}

\ifusepdf
\usepackage{hyperref}
\else\fi
\ifpsfont
\usepackage[T1]{fontenc}
\usepackage{times}
\else\fi

\iftwelvept
\else\fi

\theoremstyle{plain}
\newtheorem{Theorem}{Theorem}[section]
\newtheorem{Proposition}[Theorem]{Proposition}
\newtheorem{Lemma}[Theorem]{Lemma}
\newtheorem{Corollary}[Theorem]{Corollary}

\newtheorem{Claim}{Claim}[Theorem]

\theoremstyle{definition}

\newtheorem{Remark}[Theorem]{Remark}
\newtheorem{Example}[Theorem]{Example}

\renewcommand{\theTheorem}{\arabic{section}.\arabic{Theorem}}
\renewcommand{\theClaim}{\arabic{section}.\arabic{Theorem}.\arabic{Claim}}
\renewcommand{\theequation}{\arabic{section}.\arabic{Theorem}.\arabic{Claim}}

\def\rom{\textup}
\newcommand{\ZZ}{{\mathbb{Z}}}
\newcommand{\QQ}{{\mathbb{Q}}}
\newcommand{\RR}{{\mathbb{R}}}
\newcommand{\CC}{{\mathbb{C}}}

\newcommand{\OO}{{\mathcal{O}}}
\newcommand{\Proj}{\operatorname{Proj}}

\newcommand{\Hom}{\operatorname{Hom}}

\newcommand{\Coker}{\operatorname{Coker}}
\newcommand{\Spec}{\operatorname{Spec}}

\newcommand{\Supp}{\operatorname{Supp}}

\newcommand{\rank}{\operatorname{rk}}
\newcommand{\acherncl}{\widehat{{c}}}
\newcommand{\aPic}{\widehat{\operatorname{Pic}}}
\newcommand{\adeg}{\widehat{\operatorname{deg}}}
\newcommand{\avol}{\widehat{\operatorname{vol}}}
\newcommand{\ah}{\hat{h}^0}
\newcommand{\ahone}{\hat{h}^1}
\newcommand{\aH}{\hat{H}^0}
\newcommand{\aHone}{\hat{H}^1}
\newcommand{\achi}{\hat{\chi}}
\newcommand{\zeros}{\operatorname{div}}
\newcommand{\can}{\operatorname{can}}
\newcommand{\vol}{\operatorname{vol}}
\newcommand{\dist}{\operatorname{dist}}
\newcommand{\sub}{\operatorname{sub}}
\newcommand{\quot}{\operatorname{quot}}
\newcommand{\characteristic}{\operatorname{char}}

\newcommand{\Sym}{\operatorname{Sym}}

\newcommand{\Proof}{{\sl Proof.}\quad}
\newcommand{\QED}{{\unskip\nobreak\hfil\penalty50\quad\null\nobreak\hfil
{$\Box$}\parfillskip0pt\finalhyphendemerits0\par\medskip}}
\newcommand{\rest}[2]{\left.{#1}\right\vert_{{#2}}}

\begin{document}

\title[Continuity of volumes on arithmetic varieties]%
{Continuity of volumes on arithmetic varieties}
\author{Atsushi Moriwaki}
\address{Department of Mathematics, Faculty of Science,
Kyoto University, Kyoto, 606-8502, Japan}
\email{moriwaki@math.kyoto-u.ac.jp}
\date{\DateTime, (\Version)}
\begin{abstract}
We introduce the volume function for
$C^{\infty}$-hermitian invertible sheaves on an arithmetic variety
as an analogue of the geometric volume function.
The main result of this paper is the continuity of
the arithmetic volume function. As a consequence, we 
have the arithmetic Hilbert-Samuel formula for
a nef $C^{\infty}$-hermitian invertible sheaf.
We also give another applications, for example, 
a generalized Hodge index theorem,
an arithmetic Bogomolov-Gieseker's inequality, etc.
\end{abstract}


\maketitle



\section*{Introduction}
\renewcommand{\theTheorem}{\Alph{Theorem}}
Let $X$ be  a $d$-dimensional projective arithmetic variety
and $\aPic(X)$ the group of isomorphism classes of $C^{\infty}$-hermitian
invertible sheaves on $X$.
For $\overline{L} \in \aPic(X)$,
the volume $\avol(\overline{L})$ of $\overline{L}$ is defined by
\[
\avol(\overline{L}) = \limsup_{m \to \infty} 
\frac{\log \# \{ s \in H^0(X, mL) \mid \Vert s \Vert_{\sup} \leq 1 \}}{m^d/d!}.
\]
For example, if $\overline{L}$ is ample, then
$\avol(\overline{L}) = \adeg (\acherncl(\overline{L})^{\cdot d})$ (cf. Lemma~\ref{lem:upper:estimate:h:0}).
This is an arithmetic analogue of the volume function
for invertible sheaves on a projective variety over a field.
The geometric volume function plays a crucial role for
the birational geometry via big invertible sheaves.
In this sense, to introduce the arithmetic analogue of it is very significant.

The first important property of the volume function is 
the characterization of
a big $C^{\infty}$-hermitian invertible sheaf by the positivity of its volume
(cf. Theorem~\ref{thm:equiv:big}).
The second one is the homogeneity of the volume function, namely,
 $\avol(n \overline{L}) = n^d \avol(\overline{L})$ for all non-negative integers $n$
 (cf. Proposition~\ref{prop:avol:hom}).
By this property, it can be extended to $\aPic(X) \otimes \QQ$.
From viewpoint of arithmetic analogue,
the most important and fundamental question is the continuity of 
\[
\avol : \aPic(X) \otimes \QQ \to \RR,
\]
that is, the validity  of the formula:
\[
\lim_{\substack{\epsilon_1, \ldots, \epsilon_n \in \QQ \\ \epsilon _1\to 0, \ldots , \epsilon_n \to 0}}
\avol(\overline{L} +\epsilon_1 \overline{A}_1 + \cdots + \epsilon_n \overline{A}_n ) = \avol(\overline{L})
\]
for any $\overline{L}, \overline{A}_1, \ldots ,\overline{A}_n \in \aPic(X) \otimes \QQ$.
The main purpose of this paper is to give an affirmative answer for the above question 
(cf. Theorem~\ref{thm:cont:arithmetic:volume}).
As a consequence, we have the following arithmetic Hilbert-Samuel formula for
a nef $C^{\infty}$-hermitian invertible sheaf:

\begin{Theorem}[cf. Corollary~\ref{cor:Hilbert:Samuel:nef}]
\label{thm:Hilbert:Samuel:nef:intro}
Let $\overline{L}$ and $\overline{N}$ be $C^{\infty}$-hermitian invertible sheaves on $X$.
If $\overline{L}$ is nef, then
\[
\log \# \{ s \in H^0(X, mL+N) \mid \Vert s \Vert_{\sup} \leq 1 \}  = 
\frac{\adeg (\acherncl_1(\overline{L})^{\cdot d})}{d!}{m^d} + o(m^d) \quad(m \gg 1).
\]
In particular, $\avol(\overline{L}) = \adeg (\acherncl_1(\overline{L})^{\cdot d})$,
and $\overline{L}$ is big if and only if
$ \adeg (\acherncl_1(\overline{L})^{\cdot d}) > 0$.
\end{Theorem}

In a more general setting,
we have the following generalized Hodge index theorem:

\begin{Theorem}[cf. Theorem~\ref{thm:ineq:avol:deg}]
\label{thm:GHIT:inro}
Let $\overline{L}$ be a
$C^{\infty}$-hermitian invertible sheaf on $X$.
We assume the following:
\begin{enumerate}
\renewcommand{\labelenumi}{(\roman{enumi})}
\item
$L_{\QQ}$ is nef on $X_{\QQ}$.

\item
$c_1(\overline{L})$ is semipositive on $X(\CC)$.

\item
$L$ has moderate growth of positive even cohomologies, that is,
there are a generic resolution of singularities $\mu : Y \to X$ and
an ample invertible sheaf $A$ on $Y$ such that,
for any positive integer $n$, there is a positive integer $m_0$ such that
\[
\log \# (H^{2i}(Y, m(n\mu^*(L) + A))) = o(m^d)
\]
for all $m \geq m_0$ and for all $i > 0$.
\end{enumerate}
Then we have
an inequality $\avol(\overline{L}) \geq \adeg(\acherncl_1(\overline{L})^{\cdot d})$.
\end{Theorem}

Theorem~\ref{thm:GHIT:inro} implies that
if $L$ is nef on every geometric fiber of $X \to \Spec(\ZZ)$,
$c_1(\overline{L})$ is semipositive on $X(\CC)$, and
$\adeg(\acherncl_1(\overline{L})^{\cdot d}) > 0$, then
$\overline{L}$ is big (cf. Corollary~\ref{cor:HIT:nef}). This is a generalization of 
\cite[Corollary~(1.9)]{ZhPL}.
Moreover we can see the arithmetic Bogomolov-Gieseker's inequality
as an application  of Theorem~\ref{thm:GHIT:inro} (cf. Corollary~\ref{cor:Bogo:dim:2}).

\medskip
In the geometric case, the above Theorem~\ref{thm:Hilbert:Samuel:nef:intro} 
can be proved by using the Riemann-Roch formula and
Fujita's vanishing theorem.
In the arithmetic case, the proof in terms of
the arithmetic Riemann-Roch theorem
seems to be difficult.
Instead of it,
we prove the continuity of the volume function by direct estimates.
For this purpose,
the technical core is the following theorem,
which was inspired by Yuan's paper \cite{Yuan}.

\begin{Theorem}[cf. Theorem~\ref{thm:h:0:estimate:big}]
\label{thm:h:0:estimate:big:intro}
Let $X$ be a
projective and generically smooth arithmetic variety
of dimension $d \geq 2$.
Let $\overline{L}$ and $\overline{A}$ be $C^{\infty}$-hermitian invertible sheaves
on $X$. We assume the following:
\begin{enumerate}
\renewcommand{\labelenumi}{(\roman{enumi})}
\item
$A$ and $L +A$ are very ample over $\QQ$.

\item
The first Chern forms
$c_1(\overline{A})$ and $c_1(\overline{L} + \overline{A)}$ 
on $X(\CC)$ are positive.

\item
There is a non-zero section $s \in H^0(X, A)$ such that
the vertical component of $\zeros(s)$ is contained in
the regular locus of $X$ and that the horizontal component of
$\zeros(s)$ is smooth over $\QQ$.
\end{enumerate}
Then there are positive constants $a_0$, $C$ and $D$
depending only on $X$, $\overline{L}$ and $\overline{A}$ such that
\begin{multline*}
\log \#\{  s \in H^0(X, a L + (b-c)A) \mid \Vert s \Vert_{\sup} \leq 1 \} \\
\leq \log \#\{  s \in H^0(X, a L - cA) \mid \Vert s \Vert_{\sup} \leq 1 \} \\
+ C b a^{d-1} + D a^{d-1} \log(a)
\end{multline*}
for all integers $a, b, c$ with $a \geq b \geq c \geq 0$ and $a \geq a_0$.
\end{Theorem}

In order to explain the technical aspects of the above theorem,
let us consider it in the geometric case, namely, we assume
that $X$ is a projective smooth variety over $\CC$, and we try to estimate
\[
\Delta = h^0(X , a L + (b-c) A) -
h^0(X, a L - c A).
\]

The first elegant way:
Let us choose an infinite sequence $\{ Y_i \}_{i=1}^{\infty}$ 
of distinct smooth members of $\vert A \vert$ such that
\[
h^0(Y_i, \rest{n L + m A}{Y_i})=
h^0(Y_j, \rest{n L + m A}{Y_j})
\]
for all $i, j$ and all integers $n, m$. Then an exact sequence
\[
0 \to H^0(X, a L - c A) \to
H^0(X, a L + (b-c) A) \to 
\bigoplus_{i=1}^b H^0(Y_i, \rest{a L + (b-c) A}{Y_i})
\]
gives rise to
$\Delta
\leq b \cdot h^0(Y_1, \rest{a(L + A)}{Y_1})$.
This argument does not work in the arithmetic situation.

The second way:
In the paper \cite{Yuan}, for a fixed smooth member $Y \in \vert A \vert$,
Yuan considered an exact sequence
\[
0 \to a L + (k-1 -c)A \to a L + (k-c) A \to
\rest{a L + (k-c) A}{Y} \to 0
\]
for each $1 \leq k \leq b$, which  yields
\begin{align*}
\Delta
& \leq \sum_{k=1}^b h^0(Y, \rest{a L + (k-c) A}{Y} ) \leq 
b  \cdot h^0(Y, \rest{a (L + A)}{Y}).
\end{align*}
This second way works if we consider the arithmetic $\achi$ instead of the number of small sections.
In this way, Yuan \cite{Yuan} obtained an arithmetic analogue of a theorem of Siu.
However, if we estimate the number of small sections by using the above way,
the growth of the contribution from error terms is larger than the main term.

The third way:
An exact sequence
\[
0 \to a L - c A \to a L + (b-c) A \to
\rest{a L +(b-c) A}{bY} \to 0
\]
gives rise to
\[
\Delta \leq h^0(bY, \rest{a L + (b-c) A}{bY}).
\]
On the other hand, using exact sequences
\[
0 \to   \rest{a L + (b -c-k)A}{Y} \to \rest{a L + (b-c) A}{(k+1)Y} \to  
\rest{a L + (b-c) A}{kY} \to 0,
\]
we have
\begin{multline*}
h^0(bY, \rest{a L + (b-c) A}{bY}) \leq
\sum_{k=0}^{b-1} h^0(Y,  \rest{a L + (b -c - k) A}{Y}) \\
\leq b \cdot h^0(Y, \rest{a (L + A)}{Y}).
\end{multline*}
In the arithmetic context, the behavior of the error terms by this way is better than the second way,
so that we could get the desired estimate.
Of course, this way is very complicated because it involves non-reduced schemes.

\medskip
The paper is organized as follows:
In Section~1, we prepare several estimates of norms on complex manifolds.
In Section~2, many formulae concerning the number of small sections are
discussed. Through Section~3, we give the proof of the main technical estimate
of the number of small sections.
In Section~4, we introduce the volume function on an arithmetic variety and
consider several basic properties.
In Section~5, we prove the continuity of the volume function and
the arithmetic Hilbert-Samuel formula for a nef $C^{\infty}$-hermitian invertible sheaf.
Finally, in Section~6, we consider the generalized Hodge index theorem and
the arithmetic Bogomolov-Gieseker's inequality.

\medskip
Finally we would like to thank Prof. Mochizuki for valuable correspondences.

\renewcommand{\thesubsubsection}{\arabic{subsubsection}}

\bigskip
\renewcommand{\theequation}{CT.\arabic{subsubsection}.\arabic{Claim}}
\subsection*{Conventions and terminology}
We fix several conventions and terminology of this paper.

\subsubsection{}
\label{CT:round:up:down}
For a real number $x \in \RR$,
the {\em round-up} $\lceil x \rceil$, the {\em round-down} $\lfloor x \rfloor$ and
the {\em fractional part} $\{ x \}$ are defined by
\[
\lceil x \rceil := \min \{ k \in \ZZ \mid x \leq k \},\quad
\lfloor x \rfloor := \max \{ k \in \ZZ \mid k \leq x  \}\quad\text{and}\quad
\{ x \} = x - \lfloor x \rfloor.
\]

\subsubsection{}
\label{CT:two:norms:C:n}
For a complex vector $z = (z_1, \ldots, z_n) \in \CC^n$,
two norms $\vert z \vert$ and $\vert z \vert'$ are defined by
\[
\vert z \vert = \sqrt{\vert z_1 \vert^2 + \cdots + \vert z_n \vert^2}
\quad\text{and}\quad
\vert z \vert' = \vert z_1 \vert + \cdots + \vert z_n \vert.
\]
Note that $\vert z \vert \leq \vert z \vert' \leq \sqrt{n} \vert z \vert$ for all $z \in \CC^n$.

\subsubsection{}
\label{CT:sub:quot:norm}
Let $(V, \sigma)$ be a finite dimensional normed vector space over $\RR$.
The norm $\sigma$ is sometimes denoted by $\Vert\cdot\Vert$.
Let $f : W \to V$ be an injective homomorphism of vector spaces over $\RR$.
Then the norm $\sigma$ on $V$ yields a norm $\sigma'$ on $W$
given by $\sigma'(x) = \sigma(f(x))$. This norm $\sigma'$ is denoted by
$\sigma_{W \hookrightarrow V}$ and is called the {\em subnorm} of $\sigma$.
Let $g : V \to Q$ be a surjective homomorphism
of vector spaces over $\RR$. Then a norm $\sigma''$ on $Q$
is defined by
\[
\sigma''(y) = \inf \{ \sigma(x) \mid x \in g^{-1}(y) \}.
\]
This norm $\sigma''$ is denoted by
$\sigma_{V \twoheadrightarrow Q}$ and is called the {\em quotient norm} of $\sigma$.
Let 
\[
0 \to V' \to V \to V'' \to 0
\]
be an exact sequence of finite dimensional
vector spaces over $\RR$.
Let $\sigma'$, $\sigma$ and $\sigma''$ be norms of $V'$, $V$ and $V''$ respectively.
We say
\[
 0 \to (V', \sigma') \to (V, \sigma) \to (V'', \sigma'') \to 0
\]
is an {\em exact sequence of normed vector spaces} if
$\sigma' = \sigma_{V' \hookrightarrow V}$ and
$\sigma'' = \sigma_{V \twoheadrightarrow V''}$.
Let $V^{\vee}$ be the dual space of $V$, that is,
$V^{\vee} = \Hom_{\RR}(V, \RR)$.
The {\em dual norm} $\sigma^{\vee}$ of $V^{\vee}$ is given by
\[
\sigma^{\vee}(\phi) = \sup \{ \vert \phi(x) \vert \mid
\text{$x \in V$ and $\sigma(x) \leq 1$} \}.
\]

\subsubsection{}
\label{CT:additive:notion:invertible:sheaves}
Let $X$ be either a scheme or a complex space.
Let $L_1, \ldots, L_n$ be invertible sheaves on $X$ and
$m_1, \ldots, m_n$ integers.
In this paper, the tensor product $L_1^{\otimes m_1} \otimes \cdots \otimes L_n^{m_n}$
of invertible sheaves
is usually denoted by
\[
m_1 L_1 + \cdots + m_n L_n
\]
in the additive way like divisors.

\subsubsection{}
\label{CT:subnorm:hermitian:invertible:sheaf}
Let $X$ be a compact complex manifold and $\Omega$ a volume form on $X$.
Let $\overline{L} = (L, \vert \cdot\vert_L)$
be a $C^{\infty}$-hermitian invertible sheaf on $X$.
Then the natural {\em $L^2$-norm} $\Vert\cdot\Vert^{\overline{L}}_{L^2,\Omega}$ and
the {\em sup-norm} $\Vert\cdot\Vert^{\overline{L}}_{\sup}$
on $H^0(X, L)$ are defined by
\[
\Vert s \Vert^{\overline{L}}_{L^2,\Omega} = \left(\int_X \vert s \vert_L^2 \Omega\right)^{1/2}
\quad\text{and}\quad
\Vert s \Vert^{\overline{L}}_{\sup} = \sup \{ \vert s \vert_L(x) \mid x \in X \}
\]
for $s \in H^0(X, L)$. For simplicity,
$\Vert\cdot\Vert^{\overline{L}}_{L^2,\Omega}$ (resp. $\Vert\cdot\Vert^{\overline{L}}_{\sup}$)
is often denoted by
$\Vert\cdot\Vert^{\overline{L}}_{L^2}$ or $\Vert\cdot\Vert_{L^2}$
(resp $\Vert\cdot\Vert_{\sup}$).
For a real number $\lambda$,
a $C^{\infty}$-hermitian invertible sheaf 
$(L, \exp(-\lambda)\vert\cdot\vert_L)$
is denoted by $\overline{L}^{\lambda}$.
Let $\overline{A}$ be a positive $C^{\infty}$-hermitian invertible sheaf
on $X$.
The {\em normalized volume form $\Omega(\overline{A})$ associated with $\overline{A}$}
is given by
\[
\Omega(\overline{A}) = \frac{c_1(\overline{A})^{\wedge d}}{\int_X c_1(\overline{A})^{\wedge d}},
\]
where $c_1(\overline{A})$ is the first Chern form of $\overline{A}$ and $d = \dim X$.
Note that $\int_X \Omega(\overline{A}) = 1$.

\subsubsection{}
\label{CT:Arakelov}
A quasi-projective scheme over $\ZZ$ is called an {\em arithmetic variety}
if $X$ is an integral scheme and flat over $\ZZ$.
We say $X$ is {\em generically smooth} if $X$ is smooth over $\QQ$.
By Hironaka's resolution of singularities \cite{Hiro},
there is a projective birational morphism $\mu : X' \to X$ of
arithmetic varieties such that
$X'$ is generically smooth.
This $\mu : X' \to X$ is called a {\em generic resolution of singularities of $X$}.

\subsubsection{}
\label{CT:Arakelov:positivity}
Let $X$ be a projective arithmetic variety and
$\overline{L}$ a $C^{\infty}$-hermitian invertible sheaf on $X$.
According to \cite{MoArHt},
we define three kinds of the positivity of $\overline{L}$ as follows:

$\bullet$ {\em ample} :
$\overline{L}$ is ample if
$L$ is ample on $X$, the first Chern form $c_1(\overline{L})$ is positive on $X(\CC)$ and
$nA$ is generated by sections $s \in H^0(X, nA)$ with
$\Vert s \Vert_{\sup} < 1$ for a sufficiently large $n$.

$\bullet$ {\em nef} :
$\overline{L}$ is nef
if the first Chern form $c_1(\overline{L})$ is semipositive and
$\adeg (\rest{\overline{H}}{\Gamma}) \geq 0$ for
any $1$-dimensional closed subscheme $\Gamma$ in $X$.

$\bullet$ {\em big} :
$\overline{L}$ is big if
$L_{\QQ}$ is big on $X_{\QQ}$ and
there are a positive integer $n$ and a non-zero section 
$s$ of $H^0(X, nL)$ with $\Vert s \Vert_{\sup} < 1$.

\smallskip\noindent
By \cite[Corollary~(5.7)]{ZhPL}, if $\overline{L}$ is ample, then,
for a sufficiently large integer $n$,
$H^0(X, nL)$ has a basis $s_1, \ldots, s_N$ as a $\ZZ$-module with
$\Vert s_i \Vert_{\sup} < 1$ for all $i=1, \ldots, N$.

\subsubsection{}
\label{CT:Arakelov:order:hermitian}
Let $X$ be a projective arithmetic variety,
and let $\overline{L}$ and $\overline{M}$ be $C^{\infty}$-hermitian invertible sheaves on $X$.
We say {\em $\overline{L}$ is less than or equal to $\overline{M}$}, denoted by
$\overline{L} \leq \overline{M}$, if
there is an injective homomorphism
$\phi : L \to M$ such that
$\vert \phi_{\CC}(\cdot) \vert_M \leq \vert \cdot \vert_L$ on $X(\CC)$,
where $\vert\cdot\vert_L$ and $\vert\cdot\vert_M$ are hermitian norms
of $\overline{L}$ and $\overline{M}$ respectively. 
The following properties are easily checked (for the proof, see
Remark~\ref{rem:order:hermitian:vector:space}):
\begin{enumerate}
\renewcommand{\labelenumi}{(\arabic{enumi})}
\item 
$\overline{L} \leq \overline{M}$ if and only if $-\overline{M} \leq -\overline{L}$.

\item
If $\overline{L} \leq \overline{M}$ and $\overline{L}' \leq \overline{M}'$,
then $\overline{L} + \overline{L}'  \leq \overline{M}+ \overline{M}'$.
\end{enumerate}

\renewcommand{\theTheorem}{\arabic{section}.\arabic{subsection}.\arabic{Theorem}}
\renewcommand{\theClaim}{\arabic{section}.\arabic{subsection}.\arabic{Theorem}.\arabic{Claim}}
\renewcommand{\theequation}{\arabic{section}.\arabic{subsection}.\arabic{Theorem}.\arabic{Claim}}

\section{Several estimates of norms on complex manifolds}

\subsection{Gromov's inequality}
In this subsection, we consider Gromov's inequality and
its variants.
Let us begin with the local version of Gromov's inequality.

\begin{Lemma}[Local Gromov's inequality]
\label{lem:local:Gromov:inequality}
Let $a, b, c$ be real numbers with $a > b > c > 0$.
We set $U = \{ z \in \CC^n \mid \vert z \vert < a \}$,
$V = \{ z \in \CC^n \mid \vert z \vert < b \}$ and
$W = \{ z \in \CC^n \mid \vert z \vert < c \}$.
Let $\Omega$ be a volume form on $U$, and
let $\overline{H}_1, \ldots, \overline{H}_l$ be $C^{\infty}$-hermitian
invertible sheaves on $U$.
Let $\omega_1, \ldots, \omega_l$ be free bases of $H_1, \ldots, H_l$ over $U$
respectively.
Then there is a constant
$C$  depending only on $\overline{H}_1, \ldots, \overline{H}_l$,
$\omega_1, \ldots, \omega_l$, $\Omega$, $a$, $b$, $c$ and $n$ such that,
for any positive real number $p$, 
all non-negative integers $m_1, \ldots, m_l$
and
all $s \in H^0(U, m_1 H_1 + \cdots + m_l H_l)$,
\begin{multline*}
\max_{x \in \overline{W}} \{ \vert s \vert^p_{(m_1, \ldots, m_l)}(x) \}
\leq C(\lceil p \rceil)^{2n}(m_1 + \cdots + m_l +1)^{2n}
\left( \int_V \vert s \vert_{(m_1, \ldots, m_l)} ^p \Omega \right),
\end{multline*}
where $\vert\cdot\vert_{(m_1, \ldots, m_l)}$
is the hermitian norm of 
$m_1 \overline{H}_1 + \cdots + m_l \overline{H}_l$ and
$\lceil p \rceil$ is the round-up of $p$ 
\rom{(}cf. Conventions and terminology~\rom{\ref{CT:round:up:down}}\rom{)}.
\end{Lemma}

\Proof
Let $\vert \cdot\vert_i$ be the hermitian norm of
$\overline{H}_i$ and
$u_i = \vert \omega_i \vert_i$ on $U$.
Considering an upper bound of the partial derivatives of $u_i$ over $\overline{V}$,
we can find a positive constant $K_i$ such that
\[
\vert u_i(x) - u_i(y) \vert \leq K_i \vert x - y \vert'
\]
for all $x, y \in \overline{V}$ (for the definition of $\vert\cdot\vert'$, see
 Conventions and terminology~\ref{CT:two:norms:C:n}).
We set
\[
D = \max\left\{ \max_{ x \in \overline{V}} \left\{ \frac{K_1}{u_1(x)} \right\},
\ldots, \max_{ x \in \overline{V}} \left\{ \frac{K_l}{u_l(x)} \right\}, \frac{1}{b-c} \right\}
\quad\text{and}\quad
R = 1/D.
\]
Then, for $x_0, x \in \overline{V}$,
\begin{align*}
u_i(x) & \geq u_i(x_0) - K_i \vert x - x_0 \vert' =  
u_i(x_0) \left( 1 - \frac{K_i}{u_i(x_0)}\vert x - x_0 \vert' \right) \\
& \geq u_i(x_0) ( 1 - D \vert x - x_0 \vert').
\end{align*}
We set $B(x_0, R) = \{ x \in \CC^n \mid \vert x - x_0 \vert'  \leq R \}$.
Then $1 - D \vert x - x_0 \vert' \geq 0$ for all $x \in B(x_0, R)$.
Moreover,  if $x_0 \in \overline{W}$, then $B(x_0, R) \subseteq \overline{V}$ 
because
\[
\vert x - x_0 \vert \leq \vert x - x_0 \vert' \leq R \leq b - c.
\]
Here we claim the following: 
\begin{Claim}
\label{Claim:lem:local:Gromov:inequality:1}
For a non-negative real number $m$,
\[
\int_0^1 \cdots \int_0^1 x_1 \cdots x_n \left(1  - \frac{1}{n} (x_1 + \cdots + x_n)\right)^m dx_1 \cdots dx_n 
\geq \frac{1}{(\lceil m \rceil+1)^n(\lceil m \rceil +2)^n}.
\]
\end{Claim}

First let us consider the case where $m$ is an integer.
If $m=0$, then the assertion is obvious, so that we assume $m \geq 1$. 
Since
\begin{multline*}
\left(1  - \frac{1}{n}(x_1 + \cdots + x_n)\right)^m = \frac{1}{n^m}
\left( \sum_{i=1}^n \left(1 - x_i\right) \right)^m \\
= \frac{1}{n^m} \sum_{\substack{m_1 + \cdots + m_n = m \\ m_1 \geq 0, \ldots, m_n \geq 0}}
\frac{m!}{m_1! \cdots m_n !} (1- x_1)^{m_1} \cdots \left(1-x_n\right)^{m_n}
\end{multline*}
and
\[
\int_0^1 x(1-x)^d dx = \frac{1}{(d+1)(d+2)}
\]
for a non-negative integer $d$, the integral $I$ in the claim is equal to
\[
\frac{1}{n^m}\sum_{\substack{m_1 + \cdots + m_n = m\\ m_1 \geq 0, \ldots, m_n \geq 0}}
\frac{m!}{m_1! \cdots m_n !} \frac{1}{(m_1 +1)(m_1 + 2) \cdots (m_n +1)(m_n+2)}.
\]
Thus
\[
I \geq \frac{1}{(m+1)^n(m+2)^n n^m} 
\sum_{\substack{m_1 + \cdots + m_n = m \\ m_1 \geq 0, \ldots, m_n \geq 0}}
\frac{m!}{m_1! \cdots m_n !} = \frac{1}{(m+1)^n(m+2)^n} .
\]
If $m$ is not integer, then
\[
\left( 1  - \frac{1}{n} (x_1 + \cdots + x_n) \right)^m
\geq \left( 1  - \frac{1}{n} (x_1 + \cdots + x_n) \right)^{\lceil m \rceil}
\]
because $0 \leq 1  - \frac{1}{n} (x_1 + \cdots + x_n) \leq 1$.
Thus the claim follows.

\medskip
We choose a positive constant $e$ with
$\Omega \geq e \Omega_{can}$ on $\overline{V}$, where
\[
\Omega_{can} = 
\left(\frac{\sqrt{-1}}{2}\right)^n dz_1 \wedge d\bar{z}_1 \wedge \cdots
\wedge dz_n \wedge d\bar{z}_n.
\]
Let $s$ be an element of  
$H^0(U, m_1 H_1 + \cdots + m_l H_l)$.
Then we can find a holomorphic function $f$ over $U$ with
$s = f \omega_1^{\otimes m_1} \otimes \cdots \otimes \omega_l^{\otimes m_l}$.
We also choose $x_0 \in \overline{W}$ such that the continuous function 
$\vert s \vert_{(m_1, \ldots, m_l)}$ on 
$\overline{W}$
takes the maximum value at $x_0$.
Then
\begin{align*}
\int_V \vert s \vert_{(m_1, \ldots, m_l)}^p \Omega & \geq
e\int_{B(x_0, R)}  \vert s \vert_{(m_1, \ldots, m_l)}^p \Omega_{can} \\
& = e \int_{B(x_0, R)} \vert f \vert^p u_1^{pm_1} \cdots u_l^{pm_l} \Omega_{can} \\
&\geq e 
u_1(x_0)^{pm_1} \cdots u_l(x_0)^{pm_l} \int_{B(x_0, R)}
\vert f \vert^p (1 - D \vert x - x_0 \vert')^{m} \Omega_{can},
\end{align*}
where $m = p(m_1 + \cdots + m_l)$.
Moreover, if we set 
\[
x - x_0 = (r_1 \exp(\sqrt{-1}\theta_1), \ldots, r_n \exp(\sqrt{-1}\theta_n)),
\]
then
\begin{multline*}
 \int_{B(x_0, R)}
\vert f \vert^p (1 - D \vert x - x_0 \vert')^{m} \Omega_{can} \\
=
 \int_{\substack{r_1 + \cdots + r_n \leq R \\ r_1 \geq 0, \ldots, r_n \geq 0}} 
\left( \int_0^{2\pi} \cdots \int_0^{2\pi} \vert f \vert^p d \theta_1 \cdots d\theta_n \right)
r_1 \cdots r_n  (1 - D(r_1 + \cdots + r_n) )^{m} dr_1 \cdots dr_n.
\end{multline*}
Since $\vert f \vert^p$ is subharmonic, we have
\[
\int_0^{2\pi} \cdots \int_0^{2\pi} \vert f \vert^p d \theta_1 \cdots d\theta_n \geq (2\pi)^n 
\vert f (x_0)\vert^p.
\]
Therefore, using Claim~\ref{Claim:lem:local:Gromov:inequality:1},
\begin{multline*}
 \int_{B(x_0, R)}
\vert f \vert^p (1 - D \vert x - x_0 \vert')^{m} \Omega_{can} \\
\geq
(2\pi)^n  \vert f (x_0)\vert^p \int_{\substack{r_1 + \cdots + r_n \leq R \\ r_1 \geq 0, \ldots, r_n \geq 0}} 
r_1 \cdots r_n  (1 - D(r_1 + \cdots + r_n) )^{m} dr_1 \cdots dr_n \\
\geq
(2\pi)^n \vert f (x_0)\vert^p  \int_{[0, R/n]^n}
r_1 \cdots r_n  (1 - D(r_1 + \cdots + r_n) )^{m} dr_1 \cdots dr_n \\
\geq
\frac{(2\pi)^n \vert f (x_0)\vert^p}{(nD)^{2n}}  \frac{1}{(\lceil m \rceil +1)^n(\lceil m \rceil +2)^n}.
\end{multline*}
Gathering all calculations, if we set
$C' = e (2 \pi)^n/(nD)^{2n}$,
then
\[
\int_V \vert s \vert_{(m_1, \ldots, m_l)}^p \Omega \geq
\frac{C' \vert s(x_0) \vert^p_{(m_1, \ldots, m_l)}}{(\lceil m \rceil +1)^n(\lceil m \rceil +2)^n}.
\]
Further, since $\lceil m \rceil  \leq \lceil p \rceil(m_1 + \cdots + m_l)$,
\begin{align*}
(\lceil m \rceil +1)^n(\lceil m \rceil +2)^n & \leq ( \lceil p \rceil(m_1 + \cdots + m_l) +1)^n
( \lceil p \rceil(m_1 + \cdots + m_l) +2)^n \\
& \leq (  \lceil p \rceil(m_1 + \cdots + m_l + 1))^n
(2 \lceil p \rceil(m_1 + \cdots + m_l + 1))^n \\
& = 2^n (\lceil p \rceil)^{2n}(m_1 + \cdots + m_l + 1)^{2n}.
\end{align*}
Thus we get the lemma.
\QED

The partial results of the following corollary are found in \cite{MoBogo1} and
\cite{KMY}.

\begin{Corollary}[Gromov's inequality]
\label{cor:Gromov:ineq}
Let $M$ be an $n$-dimensional compact complex manifold,
$\Omega$ a volume form on $M$, and let
$\overline{H}_1, \ldots, \overline{H}_l$ be $C^{\infty}$-hermitian
invertible sheaves on $M$.
Then there is a constant
$C$  depending only on $\overline{H}_1, \ldots, \overline{H}_l$,
$\Omega$ and $M$ such that, for any positive real number $p$, 
all integers $m_1, \ldots, m_l$ with $m_1 \geq 0, \ldots, m_l \geq 0$,
and
all $s \in H^0(M, m_1 H_1 + \cdots + m_l H_l)$,
\begin{multline*}
\max_{x \in M} \{ \vert s \vert^p_{(m_1, \ldots, m_l)}(x) \}
\leq C (\lceil p \rceil)^{2n}(m_1 + \cdots + m_l + 1)^{2n}
\left( \int_M \vert s \vert_{(m_1, \ldots, m_l)} ^p \Omega \right).
\end{multline*}
\end{Corollary}

\Proof
We take a finite covering $\{ U_i \}_{i=1, \ldots, m}$ of $M$ with the following properties:
\begin{enumerate}
\renewcommand{\labelenumi}{(\arabic{enumi})}
\item 
$U_i$ is isomorphic to $\{ z \in \CC^n \mid \vert z \vert < 1\}$ by
using a local coordinate $z_i(x) = (z_{i1}(x), \ldots, z_{in}(x))$.
We set $V_i = \{ x \in U_i \mid \vert z_i(x) \vert < 1/2 \}$ and
$W_i = \{ x \in U_i \mid \vert z_i(x) \vert < 1/4 \}$.

\item
There are local bases $\omega_{i1}, \ldots, \omega_{il}$ of
$H_1, \ldots H_l$ over $U_i$ respectively.

\item
$\bigcup_{i=1}^m W_i = M$.
\end{enumerate}
Then our corollary follows from the local Gromov's inequality.
\QED

\begin{Corollary}
\label{cor:norm:comp:L:2:subvariety}
Let $M$ be an $n$-dimensional compact complex manifold, and
let $\overline{H}_1, \ldots, \overline{H}_l$ be $C^{\infty}$-hermitian
invertible sheaves on $M$. Let $V$ be a closed complex submanifold of $M$.
Let $\Omega_M$ and $\Omega_V$ be volume forms on $M$ and $V$
respectively. Then there is a constant $C$ such that
\[
C (m_1 + \cdots + m_l + 1)^{2n} \int_M \vert s \vert^2 \Omega_M \geq
\int_V \vert \rest{s}{V} \vert^2 \Omega_V
\]
for all non-negative integers $m_1 ,\ldots, m_l$
and all $s \in H^0(X, m_1 H_1+ \cdots + m_l H_l)$.
\end{Corollary}

\Proof
Note that
\[
\Vert s \Vert_{\sup}^2 \geq 
\Vert \rest{s}{V} \Vert_{\sup}^2 \geq \frac{\int_V \vert \rest{s}{V} \vert^2 \Omega_V}{\int_V \Omega_V}.
\]
Thus the corollary follows from Gromov's inequality.
\QED

The following lemma is due to Takuro Mochizuki,
who kindly tell us its proof.
This is a variant of Gromov's inequality.

\begin{Lemma}
\label{lem:comp:norm:X:U}
Let $X$ be an $n$-dimensional compact complex manifold and
$\omega$ a positive $(1,1)$-form on $X$.
Let
$\overline{H}_1, \ldots, \overline{H}_l$ be $C^{\infty}$-hermitian
invertible sheaves on $X$.
Then, for an open set $U$ of $X$, there are positive constants $C$, $C'$ and $D'$ such that
\[
\sup_{x \in X} \{ \vert s \vert_{(m_1, \ldots, m_l)} (x) \} \leq C^{m_1+\cdots+m_l} 
\sup_{x \in U} \{ \vert s \vert_{(m_1,\ldots,m_l)}(x) \}.
\]
and
\[
\int_X \vert s \vert_{(m_1,\ldots,m_l)}^2 \omega^{\wedge n} \leq D' \cdot {C'}^{m_1+\cdots+m_l}
\int_U \vert s \vert_{(m_1,\ldots,m_l)}^2 \omega^{\wedge n}
\]
for all non-negative  integers $m_1, \ldots, m_l$
and all 
$s \in H^0(X, m_1 H_1 + \cdots + m_l H_l)$,
where $\vert\cdot\vert_{(m_1,\ldots,m_l)}$ is 
the hermitian norm of $m_1 \overline{H}_1 + \cdots + m_l \overline{H}_l$.
\end{Lemma}

\Proof
Shrinking $U$ if necessarily, we may
identify $U$ with 
$\{ x \in \CC^n \mid \vert x \vert < 1\}$.
We set $W = \{ x \in \CC^n \mid \vert x \vert < 1/2 \}$.
In this proof, we define a Laplacian $\square_{\omega}$ by
the formula:
\[
-\frac{\sqrt{-1}}{2 \pi} \partial\bar{\partial}(g) \wedge \omega^{\wedge (n-1)} =
\square_{\omega}(g) \omega^{\wedge n}.
\]
Let $a_i$ be a $C^{\infty}$-function given by
$c_1(\overline{H}_i) \wedge \omega^{\wedge (n-1)} = a_i \omega^{\wedge n}$,
where $c_1(\overline{H}_i)$ is the first Chern form of $\overline{H}_i$.
We choose a $C^{\infty}$-function $\phi_i$ on $X$ such that
\[
\int_{X} a_i \omega^{\wedge n} = \int_{X} \phi_i \omega^{\wedge n}
\]
and that $\phi_i$ is identically zero on $X \setminus W$.
Thus we can find a $C^{\infty}$-function $F_i$ with
$\square_{\omega}(F_i) = a_i - \phi_i$.
Note that $\square_{\omega}(F_i) = a_i$ on $X \setminus W$.

Let $s \in H^0(X, m_1 H_1+ \cdots + m_l H_l)$ and we set
\[
f = \vert s \vert_{(m_1, \ldots, m_l)}^2 \exp(-(m_1F_1 + \cdots + m_l F_l)).
\]

\begin{Claim}
$\max_{ x \in X \setminus W} \{ f(x) \}
= \max_{ x \in  \partial(W)} \{ f(x) \}$.
\end{Claim}

If $f$ is a constant over $X \setminus W$, then our assertion is obvious,
so that we assume that $f$ is not a constant over $X \setminus W$.
In particular, $s \not= 0$.  Since
\[
-\frac{\sqrt{-1}}{2 \pi} \partial\bar{\partial}(\log(\vert s \vert_{(m_1, \ldots, m_l)}^2)) =
c_1(m_1 \overline{H}_1 + \cdots + m_l \overline{H}_l) = 
m_1	c_1(\overline{H}_1) + \cdots + m_l c_1(\overline{H}_l),
\]
we have $\square_{\omega}(\log(f))= 0$
on $X \setminus (W \cup \Supp(\zeros(s)))$.
Let us choose $x_0 \in X \setminus W$ such that
the $C^{\infty}$-function $f$ over $X \setminus W$
takes the maximum value at $x_0$.
Note that 
\[
 x_0  \in X \setminus (W \cup \Supp(\zeros(s))).
\]
For, if $\Supp(\zeros(s)) = \emptyset$,
then our assertion is obvious.
Otherwise, $f$ is zero
at any point of $\Supp(\zeros(s))$. 

Since $\log(f)$ is harmonic over
$X \setminus (W \cup \Supp(\zeros(s)))$,
$\log(f)$ takes the maximum value at $x_0$ and
$\log(f)$ is not a constant, we have $x_0 \in  \partial(W)$
by virtue of the maximum principle of harmonic functions.
Thus the claim follows.

\medskip
We set
\[
d_i = \min_{ x \in X \setminus W} \{ \exp(-F_i)\},\quad
D_i = \max_{x \in \partial(W)} \{ \exp(-F_i) \}\quad\text{and}\quad
C = \max_{i=1, \ldots, l} \{ D_i/d_i \}.
\]
Then
\[
d_1^{m_1} \cdots d_l^{m_l} \vert s \vert_{(m_1, \ldots, m_l)}^2  \leq f
\]
over $X \setminus W$ and
\[
f \leq D_1^{m_1} \cdots D_l^{m_l} \vert s \vert_{(m_1, \ldots, m_l)}^2
\]
over $\partial(W)$. Hence
\begin{multline*}
\max_{x \in X \setminus W} \{ \vert s \vert_{(m_1, \ldots, m_l)}^2 \} 
 \leq
C^{m_1 + \cdots + m_l} \max_{ x \in  \partial(W)} \{ \vert s \vert_{(m_1, \ldots, m_l)}^2 \} \\
\leq
C^{m_1 + \cdots + m_l} \max_{ x \in  \overline{W}} \{ \vert s \vert_{(m_1, \ldots, m_l)}^2 \}.
\end{multline*}
which implies that
\[
\max_{ x \in X} \{ \vert s \vert_{(m_1, \ldots, m_l)}^2\} \leq 
C^{m_1+\cdots+m_l} \max_{ x \in  \overline{W}} \{ \vert s \vert_{(m_1, \ldots, m_l)}^2 \}.
\]
This is the first part of the lemma.
Note that $e^x \geq x + 1$ for $x \geq 0$.
Thus, by the local Gromov's inequality (cf. Lemma~\ref{lem:local:Gromov:inequality}),
there are constants $C_1$ and $D_1$ such that
\[
 \max_{ x \in  \overline{W}} \{ \vert s \vert_{(m_1, \ldots, m_l)}^2 \} 
 \leq D_1 \cdot C_1^{m_1+\cdots+m_l} \int_U \vert s \vert_{(m_1, \ldots, m_l)}^2 \Omega
\]
for all non-negative integers $m_1, \ldots, m_l$
and all $s \in H^0(X, m_1 H_1 + \cdots + m_l H_l)$.
Therefore the second assertion follows.
\QED

\subsection{Distorsion functions}
\setcounter{Theorem}{0}
Let $X$ be an $n$-dimensional projective complex manifold and
$\Omega$ a volume form of $X$ with $\int_X \Omega = 1$.
Let $\overline{H} = (H, h)$ be a $C^{\infty}$-hermitian invertible
sheaf on $X$.
For $s, s' \in H^0(X, H)$, we set
\[
\langle s, s' \rangle_{\overline{H}, \Omega} =
\int_X h(s, s') \Omega.
\]
Let $s_1, \ldots, s_N$ be an orthonormal basis of $H^0(X, H)$
with respect to $\langle\ , \ \rangle_{\overline{H}, \Omega}$.
We define
\[
\dist(\overline{H}, \Omega)(x) =
\sum_{i=1}^N h(s_i, s_i)(x).
\]
Note that $\dist(\overline{H}, \Omega)$ does not depend on the choice
of an orthonormal basis.
In the case of $H^0(X, H) = \{ 0 \}$, $\dist(\overline{H}, \Omega)$ is defined
to be the constant function $0$.
The function $\dist(\overline{H}, \Omega)$ is called
the {\em distorsion function of $\overline{H}$ with respect to $\Omega$}.

Let $\overline{A}$ be a positive
$C^{\infty}$-hermitian invertible sheaf on $X$.
Due to Bouche \cite{Bou} and Tian \cite{Tian}, we know that
\[
\sup_{x \in X} \left\vert 
\frac{\dist(a \overline{A}, \Omega(A))(x)}{\dim H^0(a A)} - 1 \right\vert
= O(1/a)
\]
for $a \gg 1$,
where $\Omega(\overline{A})$ is the normalized volume form associated with
$\overline{A}$ (cf. Conventions and terminology~\ref{CT:subnorm:hermitian:invertible:sheaf}).
Using this result, Yuan \cite[Theorem~3.3]{Yuan} proved the following:

\begin{Theorem}
\label{thm:estimate:disfun:A:B}
Let $\overline{A} = (A, h_A)$ and $\overline{B} = (B, h_B)$ be positive $C^{\infty}$-hermitian invertible
sheaves on $X$. Then there are positive constants $C_1$ and $C_2$ such that
\[
\dist(a \overline{A} - b \overline{B}, \Omega(\overline{A}))(x)
\leq \dim H^0(a A) \left(1 + \frac{2C_1}{a} + \frac{3C_2}{b}\right)
\]
for all $x \in X$, $a \geq 1$ and $b \geq 3C_2$.
\end{Theorem}

\Proof
For reader's convenience, we reprove it here.
By Bouche-Tian's theorem,
there are constants $C_1$ and $C_2$ such that
\[
\dim H^0(aA)\left(1-\frac{C_1}{a} \right) \leq \dist(a \overline{A}, \Omega(A))(z) \leq 
\dim H^0(aA) \left( 1 + \frac{C_1}{a} \right)
\]
and
\[
\dim H^0(b B)\left( 1-\frac{C_2}{b} \right)
\leq \dist(b \overline{B}, \Omega(B))(z) \leq \dim H^0(b B) \left( 1 + \frac{C_2}{b}\right)
\]
for all $z \in X$, $a \gg 1$ and $b \gg 1$.
By taking larger $C_1$ and $C_2$ if necessarily, we may assume that
the above inequalities hold for all $z \in X$ and all $a, b \geq 1$.

Let us fix an arbitrary $x \in X$.
Let us choose an orthonormal basis of $H^0(bB)$
with respect to $\langle\ , \ \rangle_{b \overline{B}, \Omega(\overline{B})}$
such that only one section is non-zero at $x$.
We denote this section by $s(b)$.
Then
\[
h_{b\overline{B}}(s(b), s(b))(x) = \dist(b\overline{B}, \Omega(\overline{B}))(x) 
\geq \dim H^0(bB) (1 - C_2/b).
\]
On the other hand,
\[
\Vert s(b) \Vert^2_{\sup} \leq \sup_{z \in X} \dist(b\overline{B}, \Omega(B))(z) \leq
\dim H^0(bB) (1 + C_2/b).
\]
Therefore
\[
\frac{h_{b\overline{B}}(s(b), s(b))(x)}{\Vert s(b) \Vert^2_{\sup}} \geq
\frac{1 - C_2/b}{1+ C_2/b}.
\]

We choose an orthonormal basis $t_1, \ldots, t_r$ of
$H^0(aA - bB)$ with respect to
$\langle\ ,\ \rangle_{a\overline{A}- b \overline{B}, \Omega(\overline{A})}$
such that
$s(b)t_1, \ldots s(b) t_r$ is orthogonal
with respect to $\langle\ ,\ \rangle_{a\overline{A}, \Omega(\overline{A})}$
in $H^0(aA)$.
This is possible because a hermitian matrix is diagonalizable by
an unitary matrix.
Then 
\[
\{ s(b) t_i/\Vert s(b)t_i \Vert_{a\overline{A}, \Omega(\overline{A})} \}_{i=1, \ldots, r}
\]
is a part of an orthonormal basis of $H^0(aA)$.
Thus
\[
\sum_{i=1}^r \frac{h_{a \overline{A}}(s(b) t_i, s(b) t_i)(x)}
{\Vert s(b)t_i \Vert^2_{a\overline{A}, \Omega(\overline{A})}} \leq 
\dist(aA, \Omega(\overline{A}))(x) \leq \dim H^0(aA)(1 + C_1/a).
\]
On the other hand,
\[
\Vert s(b)t_i \Vert^2_{a\overline{A}, \Omega(\overline{A})} \leq
\Vert s(b) \Vert^2_{\sup} .
\]
Therefore
\begin{align*}
\frac{1-C_2/b}{1+C_2/b} \dist(a\overline{A} - b \overline{B}, \Omega(\overline{A}))(x) & \leq
\frac{h_{b\overline{B}}(s(b), s(b))(x)}{\Vert s(b) \Vert^2_{\sup}} 
\sum_{i=1}^r h_{a\overline{A}-b\overline{B}}(t_i,  t_i)(x) \\
& \leq
\sum_{i=1}^r  \frac{h_{b\overline{B}}(s(b), s(b))(x)}{\Vert s(b)t_i \Vert^2_{a\overline{A}, \Omega(\overline{A})}} 
h_{a\overline{A}-b\overline{B}}(t_i,  t_i)(x) \\
& = \sum_{i=1}^r \frac{h_{a \overline{A}}(s(b) t_i, s(b) t_i)(x)}
{\Vert s(b)t_i \Vert^2_{a\overline{A}, \Omega(\overline{A})}} \leq
\dim H^0(aA)(1 + C_1/a).
\end{align*}
Thus, if $b \geq 3C_2$, then
\[
\dist(a\overline{A} - b \overline{B}, \Omega(\overline{A}))(x) \leq
\dim H^0(aA) \frac{(1+C_1/a)(1+C_2/b)}{1-C_2/b}.
\]
It is easy to see that
\[
 \frac{(1+C_1/a)(1+C_2/b)}{1-C_2/b} = 1 + \frac{2C_1}{a} + \frac{3C_2}{b} - \frac{b-3C_2}{b-C_2} \left( \frac{C_1}{a} + \frac{C_2}{b} \right).
\]
Therefore, if $b \geq 3C_2$, then
\[
 \frac{(1+C_1/a)(1+C_2/b)}{1-C_2/b} \leq 1 + \frac{2C_1}{a} + \frac{3C_2}{b}.
\]
\QED

Let $\overline{L}$  and $\overline{A}$ be $C^{\infty}$-hermitian invertible
sheaves on a projective complex manifold $X$. 
Assume that $\overline{A}$ and $\overline{L} + \overline{A}$ are positive.
We set $\Omega = \Omega(\overline{L} + \overline{A})$.
Let $a, b, c$ be non-negative integers.
Let $s$ be a non-zero element of $H^0(bA)$
with $\Vert s \Vert_{\sup} \leq 1$.
Let $\langle\ ,\ \rangle_{a \overline{L} - c \overline{A}}$
and $\langle\ ,\ \rangle_{a \overline{L} + (b-c) \overline{A}}$ be
the natural hermitian metric of $H^0(aL - cA)$ and
$H^0(aL + (b-c)A)$ with respect to $\Omega$.
We set
\[
\begin{cases}
B_{L^2} = \{ t \in H^0(aL-cA) \mid \langle t, t \rangle_{a \overline{L} - c \overline{A}} \leq 1\} \\
B_{\sub} = \{ t \in H^0(aL-cA) \mid \langle st, st \rangle_{a \overline{L} + (b-c) \overline{A}} \leq 1\}.
\end{cases}
\]
Then we have the following corollary, which is a variant
of \cite[Proposition~3.1]{Yuan}

\begin{Corollary}
\label{cor:comp:sub:L2:ball}
There are positive constants $C_1$ and $C_2$ such that
\begin{multline*}
\log\left(\frac{\vol(B_{L^2})}{\vol(B_{\sub})}\right)
\geq \dim H^0(a(L+A)) \left(\int_X \log(\vert s \vert) \Omega \right)
\left(1 + \frac{2C_1}{a} + \frac{3C_2}{a+c}\right)
\end{multline*}
for all $a > 3C_2$ and $c \geq 0$.
\end{Corollary}

\Proof
Since $a \overline{L} - c \overline{A} = a (\overline{L} + \overline{A}) - (a + c) \overline{A}$,
by Theorem~\ref{thm:estimate:disfun:A:B},
there are positive constants $C_1$ and $C_2$ such that
\[
\dist(a \overline{L} - c \overline{A}, \Omega)(x)
\leq \dim H^0(a (L+A)) \left(1 + \frac{2C_1}{a} + \frac{2C_2}{a+c}\right)
\]
for all $x \in X$, $a \geq 1$ and $a+c \geq 3C_2$.
Note that if $a > 3 C_2$ and $c \geq 0$, then
$a+c \geq 3C_2$ and $a \geq 1$.

We choose an orthonormal basis $t_1, \ldots, t_r$ of
$H^0(aL - cA)$ with respect to
$\langle\ ,\ \rangle_{a \overline{L} - c \overline{A}}$
such that
$st_1, \ldots, st_r$ are orthogonal
with respect to $\langle\ ,\ \rangle_{a \overline{L} + (b-c) \overline{A}}$.
Then, using Jensen's inequality, for $a > 3C_2$ and $c \geq 0$,
{\allowdisplaybreaks
\begin{align*}
\log\left(\frac{\vol(B_{L^2})}{\vol(B_{\sub})}\right) & = \sum_{i=1}^r
\log \Vert s t_i \Vert_{a \overline{L} + (b-c) \overline{A}, \Omega}
= \frac{1}{2} \sum_{i=1}^r \log \int_X \vert s \vert^2 \vert t_i \vert^2 \Omega \\
& \geq  \frac{1}{2} \sum_{i=1}^r  \int_X \log(\vert s \vert^2) \vert t_i \vert^2 \Omega \\
& =  \frac{1}{2} \int_X \log(\vert s \vert^2) \dist(a \overline{L} - c \overline{A}, \Omega) \Omega \\
& \geq \dim H^0(a(L+A)) \left(\int_X \log(\vert s \vert)\Omega\right) 
\left(1 + \frac{2C_1}{a} + \frac{3C_2}{a+c}\right).
\end{align*}}
\QED

\renewcommand{\theTheorem}{\arabic{section}.\arabic{Theorem}}
\renewcommand{\theClaim}{\arabic{section}.\arabic{Theorem}.\arabic{Claim}}
\renewcommand{\theequation}{\arabic{section}.\arabic{Theorem}.\arabic{Claim}}

\section{Normed $\ZZ$-module and its invariants $\ah$, $\ahone$ and $\achi$}
Let $(M, \Vert\cdot\Vert)$ be a normed finitely generated $\ZZ$-module, namely,
$M$ is a finitely generated $\ZZ$-module and $\Vert\cdot\Vert$ is a norm
on $M_{\RR} = M \otimes_{\ZZ} \RR$.
We define $\aH(M, \Vert\cdot\Vert)$ and $\ah(M, \Vert\cdot\Vert)$ to be
\[
\aH(M, \Vert\cdot\Vert) = \{ x \in M \mid \Vert x \Vert \leq 1 \}\quad\text{and}\quad
\ah(M, \Vert\cdot\Vert) = \log \#\aH(M, \Vert\cdot\Vert).
\]
It is easy to see that
\[
\ah(M, \Vert\cdot\Vert) = \ah(M/M_{tor}, \Vert\cdot\Vert)  + 
\log \#(M_{tor}),
\]
where $M_{tor}$ is the torsion part of $M$.
We set
\[
B(M, \Vert\cdot\Vert) = \{ x \in M_{\RR} \mid \Vert x \Vert \leq 1 \}.
\]
Then $\achi(M, \Vert\cdot\Vert) $ is defined by
\[
\achi(M, \Vert\cdot\Vert) = \log \left(\frac{\vol(B(M,\Vert\cdot\Vert))}{\vol(M_{\RR}/(M/M_{tor}))}\right) +
\log \#(M_{tor}).
\]
Note that $\achi(M, \Vert\cdot\Vert)$ does not depend on the choice of
a Lebesgue measure of $M_{\RR}$ arising from  a basis of $M_{\RR}$.
Let $M^{\vee}$ be the dual of $M$, that is,
$M^{\vee} = \Hom_{\ZZ}(M, \ZZ)$.
Note that $M^{\vee}$ is torsion free.
Since $(M^{\vee})_{\RR}$ is naturally isomorphic to $(M_{\RR})^{\vee}$,
we denote $(M^{\vee})_{\RR}$ by $M^{\vee}_{\RR}$.
The norm $\Vert\cdot\Vert$ of $M_{\RR}$ yields  the dual norm $\Vert\cdot\Vert^{\vee}$
of  $M_{\RR}^{\vee}$ as follows:
for $\phi \in M_{\RR}^{\vee}$,
\[
\Vert \phi \Vert^{\vee} = \sup \{ \vert \phi(x) \vert \mid x \in B(M, \Vert\cdot\Vert) \}.
\]
Then $\aHone(M, \Vert\cdot\Vert)$ and $\ahone(M, \Vert\cdot\Vert)$ are defined by
\[
\aHone(M, \Vert\cdot\Vert) = \aH(M^{\vee}, \Vert\cdot\Vert^{\vee})\quad\text{and}\quad
\ahone(M, \Vert\cdot\Vert) = \ah(M^{\vee}, \Vert\cdot\Vert^{\vee}).
\]
Let $\Sigma = \{ e_1, \ldots. e_r \}$ be a free basis of
of $M/M_{tor}$ and let
$\langle \ , \ \rangle_{\Sigma}$ be the standard inner product of $M/M_{tor}$
in terms
of the basis $\Sigma$, that is,
\[
\langle x , y \rangle_{\Sigma} = a_1b_1 + \cdots + a_r b_r
\]
for $x = a_1 e_1,+ \cdots + a_r e_r,
y = b_1 e_1 + \cdots + b_r e_r \in M/M_{tor}$.
Then we can see
\[
\ahone(M, \Vert\cdot\Vert) = \log \# \{ x \in M/M_{tor} \mid
\text{$\vert \langle x, y \rangle_{\Sigma} \vert \leq 1$
for all $y  \in B(M, \Vert\cdot\Vert)$}\}.
\]
In the case where $M = \{ 0\}$, 
$\ah(M, \Vert\cdot\Vert)$, $\ahone(M, \Vert\cdot\Vert)$ and
$\achi(M, \Vert\cdot\Vert)$ are defined to be $0$.
The following proposition is very useful to estimate
$\ah$ of normed $\ZZ$-module.
This is essentially the results in Gillet-Soul\'{e} \cite{GSNLP}.
The following formulae are also pointed out in Yuan's paper \cite{Yuan}.

\begin{Proposition}
\label{prop:basic:property:h:0:h:1:x}
\begin{enumerate}
\renewcommand{\labelenumi}{(\arabic{enumi})}
\item For  a normed  finitely generated $\ZZ$-module $(M, \Vert\cdot\Vert)$,
\begin{multline*}
\qquad
-  \log( 6 ) \rank M \leq \ah(M, \Vert\cdot\Vert) - \ahone(M, \Vert\cdot\Vert) - \achi(M, \Vert\cdot\Vert) \\
\leq \log (3/2) \rank M + 2\log ((\rank M)!) .
\end{multline*}

\item Let $\Vert\cdot\Vert_1$ and $\Vert\cdot\Vert_2$ be two norms
of  a finitely generated $\ZZ$-module $M$ with
$\Vert\cdot\Vert_1 \leq \Vert\cdot\Vert_2$.
Then 
\[
\ah(M, \Vert\cdot\Vert_1) \geq \ah(M, \Vert\cdot\Vert_2)\quad\text{and}\quad
\ahone(M, \Vert\cdot\Vert_1) \leq \ahone(M, \Vert\cdot\Vert_2).
\]
Moreover,
\begin{multline*}
\qquad\quad
\achi(M, \Vert\cdot\Vert_2) - \achi(M, \Vert\cdot\Vert_1)
\leq  \ah(M, \Vert\cdot\Vert_2) - \ah(M, \Vert\cdot\Vert_1) \\
+ \log (9) \rank M + 2\log ((\rank M)! ).
\end{multline*}

\item
For a non-negative real number $\lambda$,
\begin{multline*}
\qquad\quad
0 \leq \ah(M, \exp(-\lambda)\Vert\cdot\Vert) - \ah(M, \Vert\cdot\Vert) \\
\leq  \lambda \rank M + \log (9) \rank M + 2\log ((\rank M)! ).
\end{multline*}

\item
Let 
\[
0 \to (M', \Vert\cdot\Vert') \overset{f}{\longrightarrow}
(M, \Vert\cdot\Vert) \overset{g}{\longrightarrow} (M'', \Vert\cdot\Vert'') \to 0
\]
be an exact sequence of  normed finitely generated $\ZZ$-modules, that is,
\[
0 \to M' \overset{f}{\longrightarrow} M
\overset{g}{\longrightarrow} 
M'' \to 0
\]
is an exact sequence of finitely generated $\ZZ$-modules and
\[
0 \to (M'_{\RR}, \Vert\cdot\Vert') \overset{f_{\RR}}{\longrightarrow}
(M_{\RR}, \Vert\cdot\Vert) \overset{g_{\RR}}{\longrightarrow} (M''_{\RR}, \Vert\cdot\Vert'') \to 0
\]
is an exact sequence of normed vector spaces over $\RR$.
Then
\begin{multline*}
\qquad\quad
\ah(M, \Vert\cdot\Vert) \leq \ah(M', \Vert\cdot\Vert') + \ah(M'', \Vert\cdot\Vert'') 
+ \log (18) \rank M' \\
+ 2\log ((\rank M')! ).
\end{multline*}


\item
If there is a basis $\{ e_1, \ldots, e_{\rank M} \}$
of $M/M_{tor}$ with $\Vert e_i \Vert \leq 1$ for all $i$,
then 
\[
\ahone(M, \Vert\cdot\Vert) \leq  \log(3) \rank M.
\]
\end{enumerate}
\end{Proposition}

\Proof
First we would like to give remarks on the paper \cite{GSNLP} due to Gillet-Soul\'{e}.
We use the same notation as in \cite{GSNLP}.
Let $K$ be a convex centrally symmetric bounded and absorbing set in $\RR^n$.
Let $K^{*}$ be the polar body of $K$, i.e.,
\[
K^{*} = \{ x \in \RR^n \mid \vert \text{$\langle x, y \rangle \vert \leq 1$ for all $y \in K$} \}.
\]
We denote the volume of $K$ by $V(K)$ and
$\#(K \cap \ZZ^n)$ by $M(K)$.
We assume an inequality 
\addtocounter{Claim}{1}
\begin{equation}
\label{eqn:prop:basic:property:h:0:h:1:x:1}
V(K)V(K^*) \geq f(n),
\end{equation}
where $f(n)$ is a constant depending only on $n$.
If we read the paper \cite{GSNLP} carefully (especially Theorem~1 and Proposition~4),
we can easily realize that the above inequality implies the following inequalities:
\addtocounter{Claim}{1}
\begin{equation}
\label{eqn:prop:basic:property:h:0:h:1:x:2}
6^{-n} \leq \frac{M(K)}{M(K^*)V(K)} \leq \frac{6^n}{f(n)}
\end{equation}
and
\addtocounter{Claim}{1}
\begin{equation}
\label{eqn:prop:basic:property:h:0:h:1:x:3}
M(K) \leq M(a K) \leq \frac{a^n M(K) 36^n}{f(n)} \quad \text{(for $a \in \RR$ with $a > 1$)}.
\end{equation}
Mahler showed \eqref{eqn:prop:basic:property:h:0:h:1:x:1} holds
for $f(n) = 4^n (n!)^{-2}$ (cf. \cite[\S14, Theorem~4]{GLGN}).
Bourgain and Milman \cite{BM} also proved \eqref{eqn:prop:basic:property:h:0:h:1:x:1}
for $f(n) = c^n V_n$,
where $c$ is an absolute constant and $V_n$ is the volume of the unit sphere
in $\RR^n$.
Here we uses Mahler's result for its simplicity.

\medskip
(1) 
Since
\begin{multline*}
\ah(M, \Vert\cdot\Vert) - \ahone(M, \Vert\cdot\Vert) - \achi(M, \Vert\cdot\Vert)  \\
= \ah(M/M_{tor}, \Vert\cdot\Vert) - \ahone(M/M_{tor}, \Vert\cdot\Vert) - \achi(M/M_{tor}, \Vert\cdot\Vert),
\end{multline*}
we may assume that $M$ is torsion free.
Thus (1) is a consequence of  \eqref{eqn:prop:basic:property:h:0:h:1:x:2}.

\medskip
(2) The inequalities $\ah(M, \Vert\cdot\Vert_1) \geq \ah(M, \Vert\cdot\Vert_2)$ and
$\ahone(M, \Vert\cdot\Vert_1) \leq \ahone(M, \Vert\cdot\Vert_2)$ are obvious
by their definitions. The third inequality is a consequence of (1).

\medskip
(3) Since
\begin{multline*}
\ah(M, \exp(-\lambda)\Vert\cdot\Vert) - \ah(M, \Vert\cdot\Vert) \\
= \ah(M/M_{tor}, \exp(-\lambda)\Vert\cdot\Vert) - \ah(M/M_{tor}, \Vert\cdot\Vert),
\end{multline*}
we may assume that $M$ is torsion free.
Thus
it follows from \eqref{eqn:prop:basic:property:h:0:h:1:x:3}.

\medskip
(4) We may assume $M'$ is a sub-module of $M$.
Let us choose $x_1, \ldots, x_l \in M$ with the following properties:
\begin{enumerate}
\renewcommand{\labelenumi}{(\roman{enumi})}
\item $\Vert x_i \Vert \leq 1$ for all $i$.

\item
$g(x_i) \not= g(x_j)$ for all $i\not=j$.

\item
For any $x \in M$ with $\Vert x \Vert \leq 1$,
there is $x_i$ such that $g(x) = g(x_i)$.
\end{enumerate}
By using (i) and (ii), for any $x \in M$ with $\Vert x \Vert \leq 1$,
there is a unique $x_i$ with $g(x) = g(x_i)$.
Moreover $x - x_i \in M'$ and $\Vert x - x_i\Vert \leq 2$.
On the other hand,
$\log(l) \leq \ah(M'',\Vert\cdot\Vert'')$ because $\Vert g(x_i) \Vert'' \leq 1$ for all $i$.
Therefore,
\[
\ah(M, \Vert\cdot\Vert) \leq \ah(M'',\Vert\cdot\Vert'') + \log \#\{ x' \in M' \mid \Vert x' \Vert \leq 2 \}
\]
Hence (4) follows from (3).


\medskip
(5) Let $\langle\ , \ \rangle$ be an inner product of $M/M_{tor}$
with respect to the basis $\{ e_1, \ldots, e_{\rank M} \}$.
Then, for $x = a_1 e_1 + \cdots + a_{\rank M} e_{\rank M}$,
if $\vert \langle x, e_i \rangle \vert \leq 1$ for all $i$,
then $\vert a_i \vert \leq 1$ for all $i$.
Thus (5) follows.
\QED

\begin{Remark}
\label{rem:prop:basic:property:h:0:h:1:x}
Note that 
\[
\begin{cases}
\text{$(x+1)\log(x+1) \geq x$ for all $x \geq 0$}, \\
\text{$\log(n!) \leq (n+1)\log(n+1)$ for all non-negative integer $n$}.
\end{cases}
\]
Therefore,
we have simpler inequalities for each case of Proposition~\ref{prop:basic:property:h:0:h:1:x}
as follows.
The inequalities
\eqref{eqn:rem:prop:basic:property:h:0:h:1:x:1},
\eqref{eqn:rem:prop:basic:property:h:0:h:1:x:2},
\eqref{eqn:rem:prop:basic:property:h:0:h:1:x:3} and
\eqref{eqn:rem:prop:basic:property:h:0:h:1:x:4}
are simpler versions of the corresponding inequalities in (1), (2), (3) and (4)
of Proposition~\ref{prop:basic:property:h:0:h:1:x}
respectively.

\addtocounter{Claim}{1}
\begin{multline}
\label{eqn:rem:prop:basic:property:h:0:h:1:x:1}
\left\vert \ah(M, \Vert\cdot\Vert) - \ahone(M, \Vert\cdot\Vert) - \achi(M, \Vert\cdot\Vert) 
\right\vert \\
\leq (\log(3/2) + 2)\left(\rank M + 1\right)\log\left(\rank M+ 1 \right).
\end{multline}

\addtocounter{Claim}{1}
\begin{multline}
\label{eqn:rem:prop:basic:property:h:0:h:1:x:2}
\achi(M, \Vert\cdot\Vert_2) - \achi(M, \Vert\cdot\Vert_1)
\leq  \ah(M, \Vert\cdot\Vert_2) - \ah(M, \Vert\cdot\Vert_1) \\
+ 
(\log(9) + 2) \left(\rank M + 1\right)\log\left(\rank M + 1 \right).
\end{multline}

\addtocounter{Claim}{1}
\begin{multline}
\label{eqn:rem:prop:basic:property:h:0:h:1:x:3}
0 \leq \ah(M, \exp(-\lambda)\Vert\cdot\Vert) - \ah(M, \Vert\cdot\Vert) \\
\leq  \lambda \rank M  + 
(\log(9) + 2) \left(\rank M + 1\right)\log\left(\rank M + 1 \right).
\end{multline}

\addtocounter{Claim}{1}
\begin{multline}
\label{eqn:rem:prop:basic:property:h:0:h:1:x:4}
\ah(M, \Vert\cdot\Vert) \leq \ah(M', \Vert\cdot\Vert') + \ah(M'', \Vert\cdot\Vert'') \\
+ (\log(18) + 2)\left(\rank M' + 1\right)\log\left(\rank M' + 1 \right).
\end{multline}

\end{Remark}

\section{Approximation of the number of small sections}

In this section,
we prove the main technical tool of this paper.
First we consider the following three lemmas.
The first one is 
an upper estimate of the number of small sections.

\begin{Lemma}
\label{lem:upper:estimate:h:0}
Let $X$ be a
projective arithmetic variety
of dimension $d$, and let
$\overline{L}$ and $\overline{N}$ be $C^{\infty}$-hermitian invertible sheaves on $X$.
Then we have the following:
\begin{enumerate}
\renewcommand{\labelenumi}{(\arabic{enumi})}
\item
If $\overline{L}$ is ample, then
\[
\ah\left(H^0(X, mL+N), \Vert\cdot\Vert_{\sup}^{m\overline{L}+\overline{N}}\right) =
\frac{\adeg(\acherncl_1(\overline{L})^{\cdot d})}{d!} m^d + o(m^d)
\]
for $m \gg 1$.

\item
In general, there is a constant $C$ with
\[
\ah\left(H^0(X, mL+N), \Vert\cdot\Vert_{\sup}^{m\overline{L}+\overline{N}}\right) \leq C m^d
\]
for all $m \geq 1$.

\item
Let $\mu : Y \to X$ be a generic resolution of singularities of $X$.
Let $\Omega$ be a volume form on $Y(\CC)$.
An $L^2$-norm of $H^0(X, mL + N)$ is given in the following way:
for $t \in H^0(X, mL + N)$,
\[
\Vert t \Vert_{L^2,\Omega}^{m \overline{L} + \overline{N}} := \left(
\int_{Y(\CC)} \mu^*(\vert t \vert^2_{m\overline{L} + \overline{N}}) \Omega\right)^{1/2},
\]
where $\vert\cdot\vert_{m\overline{L}+\overline{N}}$ is the hermitian norm of
$m \overline{L} + \overline{N}$.
Then there is a constant $C$ with
\[
\ah\left(H^0(X, mL+N), \Vert\cdot\Vert_{L^2, \Omega}^{m\overline{L}+\overline{N}}\right) \leq C m^d
\]
for all $m \geq 1$.
\end{enumerate}
\end{Lemma}

\Proof
(1) It is well-known that
\[
\achi\left(H^0(X, mL+N), \Vert\cdot\Vert_{\sup}^{m\overline{L}+\overline{N}}\right) =
\frac{\adeg(\acherncl_1(\overline{L})^{\cdot d})}{d!} m^d + o(m^d)
\]
for $m \gg 1$ (cf. \cite{GSNLP}, \cite{AB} and \cite[Theorem~(1.4)]{ZhPL}).
Thus, by (2) of Proposition~\ref{prop:basic:property:h:0:h:1:x},
\begin{multline*}
\ah\left(H^0(X, mL+N), \Vert\cdot\Vert_{\sup}^{m\overline{L}+\overline{N}}\right) -
\ahone\left(H^0(X, mL+N), \Vert\cdot\Vert_{\sup}^{m\overline{L}+\overline{N}}\right) \\
=
\frac{\adeg(\acherncl_1(\overline{L})^{\cdot d})}{d!} m^d + o(m^d).
\end{multline*}
Since $\overline{L}$ is ample, 
by \cite[Theorem~(4.2)]{ZhPL},
$H^0(X, mL+N)$ is generated by sections $t$
with $\Vert t \Vert_{\sup} < 1$.
Thus,
by (5) of Proposition~\ref{prop:basic:property:h:0:h:1:x},
\[
\ahone\left(H^0(X, mL), \Vert\cdot\Vert_{\sup}^{m\overline{L}}\right) = o(m^d).
\]
Hence we get (1).

\medskip
(2)
Let $\overline{A}$ be an ample $C^{\infty}$-hermitian invertible sheaf on $X$.
Then there are a positive integer $n$ and a non-zero section
$t$ of $H^0(nA - L)$ with $\Vert s\Vert_{\sup} \leq 1$.
Let $\phi : L \to nA$ be an injective homomorphism
given by $\phi(t) = s \otimes t$.
Then since $\vert s \otimes t \vert = \vert s \vert \vert t \vert \leq \vert t \vert$,
$\phi$ yields $\overline{L} \leq n\overline{A}$ (cf. Conventions and terminology~\ref{CT:Arakelov:order:hermitian}).
Therefore $m\overline{L} \leq mn\overline{A}$ for all $m \geq 1$.
Thus (2) follows from (1).

\medskip
(3)
By using Gromov's inequality on $Y(\CC)$,
there is a constant $C_1$ such that
\[
\Vert\cdot\Vert^{m \overline{L}+\overline{N}}_{\sup} \leq
C_1 m^{d-1} \Vert\cdot\Vert^{m \overline{L}+\overline{N}}_{L^2,\Omega}
\]
for $m \gg 1$.
Thus
\[
\ah\left(H^0(X, mL+N), \Vert\cdot\Vert^{m \overline{L}+\overline{N}}_{\sup} \right)
\geq \ah\left(H^0(X, mL+N), C_1 m^{d-1} \Vert\cdot\Vert^{m \overline{L}+\overline{N}}_{L^2,\Omega} \right).
\]
Moreover, by (3) of Proposition~\ref{prop:basic:property:h:0:h:1:x},
\begin{multline*}
\ah\left(H^0(X, mL+N), C_1 m^{d-1} \Vert\cdot\Vert^{m \overline{L}+\overline{N}}_{L^2,\Omega} \right) \\
=
 \ah\left(H^0(X, mL+N), \Vert\cdot\Vert^{m \overline{L}+\overline{N}}_{L^2,\Omega} \right) + o(m^{d}).
\end{multline*}
Therefore we get (3).
\QED

Next we consider formulae concerning subnorms and quotient norms
(cf. Conventions and terminology~\ref{CT:sub:quot:norm}).

\begin{Lemma}
\label{lem:sub:quot:vs:quot:sub}
\begin{enumerate}
\renewcommand{\labelenumi}{(\arabic{enumi})}
\item
Let $f : V \to W$ and $g : W \to U$ be surjective homomorphisms of finite
dimensional vector spaces over $\RR$. For a norm $\sigma$ of $V$,
$(\sigma_{V \twoheadrightarrow W})_{W \twoheadrightarrow U} = \sigma_{V \twoheadrightarrow U}$
as norms of $U$.

\item
Let
\[
\begin{CD}
W @>{f}>> V \\
@V{g'}VV @VV{g}V \\
P @>{f'}>> Q
\end{CD}
\]
be a commutative diagram of finite dimensional vector spaces over $\RR$
such that $f$ and $f'$ are injective and that $g$ and $g'$ are surjective.
Let $\sigma$ be a norm of $V$.
Then
\[
(\sigma_{W \hookrightarrow V})_{W \twoheadrightarrow P}
\geq (\sigma_{V  \twoheadrightarrow Q})_{P \hookrightarrow Q}
\]
as norms of $P$.
Moreover, if $\ker(g) \subseteq f(W)$, then
\[
(\sigma_{W \hookrightarrow V})_{W \twoheadrightarrow P}
= (\sigma_{V  \twoheadrightarrow Q})_{P \hookrightarrow Q}.
\]
\end{enumerate}
\end{Lemma}

\Proof
(1) Let us fix $u \in U$.
For $v \in (g \circ f)^{-1}(u)$, 
\[
\sigma(v) \geq \sigma_{V \twoheadrightarrow W}(f(v)) \geq 
(\sigma_{V \twoheadrightarrow W})_{W \twoheadrightarrow U}(u)
\]
Therefore $\sigma_{V \twoheadrightarrow U}(u) \geq 
(\sigma_{V \twoheadrightarrow W})_{W \twoheadrightarrow U}(u)$.

Pick up $v_0 \in (g \circ f)^{-1}(u)$ with $\sigma_{V \twoheadrightarrow U}(u) = \sigma(v_0)$.
Then, for any $w \in g^{-1}(u)$, $\sigma(v_0) \leq \sigma_{V \twoheadrightarrow W}(w)$
because $f^{-1}(w) \subseteq (g \circ f)^{-1}(u)$. Hence
\[
\sigma_{V \twoheadrightarrow U}(u) =\sigma(v_0) \leq 
(\sigma_{V \twoheadrightarrow W})_{W \twoheadrightarrow U}(u).
\]

\medskip
(2)
Since $f(\ker(g')) = f(W) \cap \ker(g)$,
for $w \in W$, 
\[
\begin{cases}
(\sigma_{W \hookrightarrow V})_{W \twoheadrightarrow P}(g'(w))
= \inf \{ \sigma(x) \mid x \in f(w) + f(W) \cap \ker(g) \} \\
(\sigma_{V  \twoheadrightarrow Q})_{g(W) \hookrightarrow Q}(g'(w)) =
\inf \{ \sigma(x) \mid x \in f(w) + \ker(g) \}.
\end{cases}
\]
Thus $(\sigma_{W \hookrightarrow V})_{W \twoheadrightarrow P}(g(w))
\geq (\sigma_{V  \twoheadrightarrow Q})_{P \hookrightarrow Q}(g(w))$.
Moreover, if $\ker(g) \subseteq f(W)$ (or, equivalently $f(W) \cap \ker(g) = \ker(g)$), then
$(\sigma_{W \hookrightarrow V})_{W \twoheadrightarrow P}
= (\sigma_{V  \twoheadrightarrow Q})_{P \hookrightarrow Q}$.
\QED

The following lemma is needed to find a good $\overline{A}$ in the proof
of Theorem~\ref{thm:h:0:estimate:big}.

\begin{Lemma}
\label{lem:reduction:assertion:theorem}
Let $X$ be a
projective and generically smooth arithmetic variety
of dimension $d$, and let $\Omega$ be a volume form on $X(\CC)$.
Let $\overline{L}$ and $\overline{A}$ be $C^{\infty}$-hermitian invertible sheaves
on $X$.  
Let us consider the following assertion $\Sigma(X, \overline{L}, \overline{A})$\rom{:}
\begin{quote}
There are positive constants $a_0$, $C$ and $D$
depending only on $X$, $\overline{L}$ and $\overline{A}$ such that
\begin{multline*}
\qquad\quad
\ah\left(H^0(a L + (b-c)A),
\Vert\cdot\Vert^{a \overline{L} + (b-c) \overline{A}}_{L^2,\Omega}\right) \\
\leq \ah\left(H^0(a L - cA), 
\Vert\cdot\Vert^{a \overline{L} - c \overline{A}}_{L^2,\Omega}\right)
+ C b a^{d-1} + D a^{d-1} \log(a)
\end{multline*}
for all integers $a, b, c$ with $a \geq b \geq c \geq 0$ and $a \geq a_0$.
\end{quote}
Then we have the following:
\begin{enumerate}
\renewcommand{\labelenumi}{(\arabic{enumi})}
\item
Let $\overline{A}'$ be another $C^{\infty}$-hermitian invertible sheaf on $X$
with $\overline{A}' \leq \overline{A}$
\rom{(}cf. Conventions and terminology~\rom{\ref{CT:Arakelov:order:hermitian}}\rom{)}.
If $\Sigma(X, \overline{L}, \overline{A})$ holds, then
so does  $\Sigma(X, \overline{L}, \overline{A}')$.

\item
We assume that  $\rank H^0(X, A) \not= 0$.
Let $\vert\cdot\vert_A$ be the hermitian norm of $\overline{A}$.
If $\Sigma(X, \overline{L}, \overline{A})$ holds, then
so does $\Sigma(X, \overline{L}, (A, \exp(-\lambda)\vert\cdot\vert_A))$ for
all $\lambda \geq 0$.

\item
We assume that $\rank H^0(X, A) \not= 0$.
Let $\overline{A}'$ be another $C^{\infty}$-hermitian invertible sheaf on $X$
such that $A'$ is isomorphic to $A$ over $\QQ$. Then
$\Sigma(X, \overline{L}, \overline{A})$ holds if and only if
so does $\Sigma(X, \overline{L}, \overline{A}')$.
\end{enumerate}
\end{Lemma}

\Proof
(1)
Since 
\[
\overline{L} + (b-c) \overline{A}' \leq \overline{L} + (b-c) \overline{A}\quad\text{and}\quad
a \overline{L} - c \overline{A} \leq a \overline{L} - c \overline{A}',
\]
(1) follows.

\medskip
(2) We set $\overline{A}' = (A, \exp(-\lambda)\vert\cdot\vert_A)$.
Let us fix constants $C_1$ and $C_2$
such that
\[
\rank H^0(a(L+A)) \leq C_1 a^{d-1}
\]
for all $a \geq 1$
and that
\[
(\log(18) + 2) \left(\rank H^0(a(L+A))+ 1\right)\log\left(\rank H^0(a(L+A))+ 1 \right)
\leq C_2 a^{d-1} \log(a)
\]
for all $a \geq 2$.
It is easy to see that
\[
\begin{cases}
\Vert\cdot\Vert^{a \overline{L} + (b-c) \overline{A}'}_{L^2,\Omega}
= \exp(-(b-c)\lambda)\Vert\cdot\Vert^{a \overline{L} + (b-c) \overline{A}}_{L^2,\Omega}, \\
\Vert\cdot\Vert^{a \overline{L} - c \overline{A}'}_{L^2,\Omega} =
\exp(c\lambda)\Vert\cdot\Vert^{a \overline{L} - c \overline{A}}_{L^2,\Omega}.
\end{cases}
\]
Since
\[
\rank H^0(aL - cA) \leq \rank H^0(aL + (b-c)A) \leq \rank H^0(a(L+A)),
\]
using \eqref{eqn:rem:prop:basic:property:h:0:h:1:x:3},
we have
\begin{multline*}
0 \leq \ah\left(H^0(a L + (b-c)A),
\Vert\cdot\Vert^{a \overline{L} + (b-c) \overline{A}'}_{L^2,\Omega}\right) \\
 - \ah\left(H^0(a L + (b-c)A),
\Vert\cdot\Vert^{a \overline{L} + (b-c) \overline{A}}_{L^2, \Omega}\right) \\
\leq C_1 \lambda (b-c) a^{d-1} + C_2 a^{d-1} \log(a)
\leq C_1 \lambda b a^{d-1} + C_2 a^{d-1} \log(a)
\end{multline*}
and
\begin{multline*}
0 \leq \ah\left(H^0(a L - cA), 
\Vert\cdot\Vert^{a \overline{L} - c \overline{A}}_{L^2,\Omega}\right)
- \ah\left(H^0(a L - cA),
\Vert\cdot\Vert^{a \overline{L} - c \overline{A'}}_{L^2,\Omega}\right) \\
\leq C_1 \lambda c a^{d-1} + C_2 a^{d-1}\log(a)
\leq C_1 \lambda b a^{d-1} + C_2 a^{d-1}\log(a).
\end{multline*}
Thus we have
\begin{multline*}
\ah\left(H^0(a L + (b-c)A),
\Vert\cdot\Vert^{a \overline{L} + (b-c) \overline{A}'}_{L^2,\Omega}\right) \leq
\ah\left(H^0(a L - cA), 
\Vert\cdot\Vert^{a \overline{L} - c \overline{A}'}_{L^2,\Omega}\right) \\
+ (C + 2 \lambda C_1)b a^{d-1} + (D + 2 C_2) a^{d-1} \log(a).
\end{multline*}
for all integers $a, b, c$ with $a \geq b \geq c \geq 0$ and $a \geq a_0$.

\medskip
(3)
It is sufficient to show that if $\Sigma(X, \overline{L}, \overline{A})$ holds,
then so does $\Sigma(X, \overline{L}, \overline{A}')$.
Since $A'$ is isomorphic to $A$ over $\QQ$,
there is a Cartier divisor $F$ such that
$A' \otimes  \OO_X(F) \simeq A$ and $\Supp(F)$ is vertical.
Thus there is a positive integer $N$ such that
$\OO_X \cdot N \subseteq \OO_X(F)$.
Hence we have a natural injective homomorphism
$\alpha : A' \cdot N \to A$.
Let $\vert\cdot\vert$ and $\vert\cdot\vert'$ be $C^{\infty}$-hermitian norms of
$\overline{A}$ and $\overline{A}'$.
Then $(A' \cdot N, \ \vert\cdot\vert')$ is a $C^{\infty}$-hermitian invertible sheaf on $X$.
Since $\alpha : A' \cdot N \to A$ is isomorphism over $\QQ$,
there is a positive number $\lambda$ such that
$\vert \alpha_{\CC}(\cdot) \vert \leq \exp(\lambda) \vert \cdot \vert'$.
Then $(A' \cdot N, \ \exp(\lambda)\vert\cdot\vert') \leq (A, \vert\cdot\vert)$.
Hence, by (1), 
$\Sigma(X, \overline{L}, (A' \cdot N, \ \exp(\lambda)\vert\cdot\vert'))$ holds.
Note that the homomorphism $A' \to A' \cdot N$ given by
$a \mapsto a \cdot N$ yields to an isometry
$(A', \ N\exp(\lambda)\vert\cdot\vert') \to (A' \cdot N, \ \exp(\lambda)\vert\cdot\vert')$.
Therefore $\Sigma(X, \overline{L}, (A', \ N\exp(\lambda)\vert\cdot\vert'))$ holds, so that
so does $\Sigma(X, \overline{L}, (A', \vert\cdot\vert'))$ by (2).
\QED

Let $X$ be a compact complex manifold, and
let $\overline{L} = (L, \vert\cdot\vert_L)$ and
$\overline{M} = (M, \vert \cdot\vert_M)$
be $C^{\infty}$-hermitian invertible sheaves on $X$.
Let $t$ be a non-zero global section of $H^0(X, M)$.
We denote by $\Vert\cdot\Vert^{\overline{L}, L - M}_{L^2, t, \sub}$
the subnorm of $H^0(X, L - M)$ induced by the natural injective homomorphism
$H^0(X, L - M) \overset{\otimes t}{\longrightarrow} H^0(X, L)$
and the $L^2$-norm of $\Vert\cdot\Vert^{\overline{L}}_{L^2}$ of 
$H^0(X, L)$
for a fixed volume form on $X$.
For simplicity,  $\Vert\cdot\Vert^{\overline{L}, L - M}_{L^2, t, \sub}$ is often
denoted by  $\Vert\cdot\Vert^{\overline{L}}_{L^2, t, \sub}$.

\medskip
The following theorem is the technical core of this paper.
The similar result for an arithmetic curve will be treated in
Proposition~\ref{prop:0:estimate:big:curve}.

\begin{Theorem}
\label{thm:h:0:estimate:big}
Let $X$ be a
projective and generically smooth arithmetic variety
of dimension $d \geq 2$.
Let $\overline{L}$ and $\overline{A}$ be $C^{\infty}$-hermitian invertible sheaves
on $X$. We assume the following:
\begin{enumerate}
\renewcommand{\labelenumi}{(\roman{enumi})}
\item
$A$ and $L +A$ are very ample over $\QQ$.

\item
The first Chern forms
$c_1(\overline{A})$ and $c_1(\overline{L} + \overline{A)}$ 
on $X(\CC)$ are positive.

\item
There is a non-zero section $s \in H^0(X, A)$ such that
the vertical component of $\zeros(s)$ is contained in
the regular locus of $X$ and that the horizontal component of
$\zeros(s)$ is smooth over $\QQ$.
\end{enumerate}
Then there are positive constants $a_0$, $C$ and $D$
depending only on $X$, $\overline{L}$ and $\overline{A}$ such that
\begin{multline*}
\ah\left(H^0(a L + (b-c)A),
\Vert\cdot\Vert^{a \overline{L} + (b-c) \overline{A}}_{L^2}\right) \leq
\ah\left(H^0(a L - cA), 
\Vert\cdot\Vert^{a \overline{L} - c \overline{A}}_{L^2}\right) \\
+ C b a^{d-1} + D a^{d-1} \log(a)
\end{multline*}
for all integers $a, b, c$ with $a \geq b \geq c \geq 0$ and $a \geq a_0$,
where the volume form $\Omega$ to define $L^2$-norms is $\Omega(\overline{L} + \overline{A})$
\rom{(}cf. Conventions and terminology~\rom{\ref{CT:subnorm:hermitian:invertible:sheaf}}\rom{)}.
Moreover the sup-version of the above estimate holds as follows:
there are positive constants $a'_0$, $C'$ and $D'$
depending only on $X$, $\overline{L}$ and $\overline{A}$ such that
\begin{multline*}
\ah\left(H^0(a L + (b-c)A),
\Vert\cdot\Vert^{a \overline{L} + (b-c) \overline{A}}_{\sup}\right) \leq
\ah\left(H^0(a L - cA), 
\Vert\cdot\Vert^{a \overline{L} - c \overline{A}}_{\sup}\right) \\
+ C' b a^{d-1} + D' a^{d-1} \log(a)
\end{multline*}
for all integers $a, b, c$ with $a \geq b \geq c \geq 0$ and $a \geq a'_0$
\end{Theorem}

\Proof
First let us fix constants $C_1$ and $C_2$
such that
\[
\rank H^0(a(L+A)) \leq C_1 a^{d-1}
\]
for all $a \geq 1$
and that
\[
(\log(18) + 2) \left(\rank H^0(a(L+A))+ 1\right)\log\left(\rank H^0(a(L+A))+ 1 \right)
\leq C_2 a^{d-1} \log(a)
\]
for all $a \geq 2$.

Let $\vert\cdot\vert_A$ be the $C^{\infty}$-hermitian norm of $\overline{A}$.
As in Conventions and terminology~\ref{CT:subnorm:hermitian:invertible:sheaf},
for $\lambda \in \RR$, we set 
\[
\overline{A}^{\lambda} = (A, \exp(-\lambda)\vert\cdot\vert_A).
\]

First we claim the following:

\begin{Claim}
\label{claim:thm:h:0:estimate:big:2}
We may assume that there is a non-zero section $s \in H^0(X, A)$ such that
$\Vert s \Vert_{\sup} \leq 1$, $\zeros(s)$ is smooth over $\QQ$ and that
$\zeros(s)$ has no vertical components.
We may further assume that there are a positive integer $n$ and
a non-zero section $t$ of $H^0(X, nA - L)$ such that
$\Vert t \Vert_{\sup} \leq 1$ and $t$ is not zero on $\zeros(s)$.
\end{Claim}

By our assumption (iii), 
there is a non-zero section $s \in H^0(X, A)$ such that
the vertical component of $\zeros(s)$ is contained in
the regular locus of $X$ and that the horizontal component of
$\zeros(s)$ is smooth over $\QQ$.
Let $Y$ and $F$ be the horizontal component of $\zeros(s)$ and
the vertical component of $\zeros(s)$ respectively.
Note that $Y$ and $F$ are effective Cartier divisors because
$F$ is contained in the regular locus of $X$.
We define a $C^{\infty}$-hermitian invertible sheaf
$\overline{A}_1$ by the equation 
\[
\overline{A} = \overline{A}_1\otimes (\OO_X(F), \vert\cdot\vert_{\can}).
\]
Then there is a non-zero section $s_1 \in H^0(X, A_1)$ such that
$s = s_1 \otimes 1_F$ and $\zeros(s_1) = Y$,
where $1_F$ is the canonical section of $\OO_X(F)$.
Let $\lambda$ be a non-negative real number with
$\exp(-\lambda) \Vert s_1 \Vert_{\sup} \leq 1$.
Then, by (3) of Lemma~\ref{lem:reduction:assertion:theorem}, 
if the assertion holds for $\overline{L}$ and $\overline{A}_1^{\lambda}$,
then so does for $\overline{L}$ and $\overline{A}$.

Moreover,
since $A$ is very ample over $\QQ$ and $\zeros(s)$
has no vertical components,
there are a positive integer $n$ and a non-zero section $t$
of $H^0(nA - L)$ such that $t$ is not zero on $\zeros(s)$.
Let $\lambda'$ be a non-negative real number
with $\exp(-\lambda'/n)\Vert t \Vert_{\sup} \leq 1$.
Then $t$ is a small section of a $C^{\infty}$-hermitian invertible sheaf
$n \overline{A}^{\lambda'} - \overline{L}$.
Thus, by (1) of Lemma~\ref{lem:reduction:assertion:theorem},
the claim follows.

\medskip
For a coherent sheaf ${\mathcal F}$ on $X$ and a subscheme $Z$ of $X$,
the image $H^i(X, {\mathcal F}) \to H^i(Z, \rest{{\mathcal F}}{Z})$ is denoted by 
$I^i(Z, \rest{{\mathcal F}}{Z})$.

If $b = 0$, then $c=0$.
Thus, in this case, the assertion is obvious, so that
we may assume $b \geq 1$.
As in Claim~\ref{claim:thm:h:0:estimate:big:2},
let $s$ be a non-zero section $H^0(X, A)$ such that
$\Vert s \Vert_{\sup} \leq 1$, $Y := \zeros(s)$ is smooth over $\QQ$ and that
$Y$ has no vertical components.
Let us choose positive numbers $C_3$ and $C_4$
such that
\[
\rank H^0(Y, \rest{a(L+A)}{Y}) \leq C_3 a^{d-2}
\]
for all $a \geq 1$
and that
\begin{multline*}
(\log(18) + 2) \left(\rank H^0(Y, \rest{a(L+A)}{Y}) + 1\right)
\log\left(\rank H^0(Y, \rest{a(L+A)}{Y}) + 1 \right) \\
\leq C_4 a^{d-2} \log(a)
\end{multline*}
for all $a \geq 2$

Let $\Vert\cdot\Vert_{L^2,\quot}^{a\overline{L}+(b-c)\overline{A}}$ 
be the quotient norm of $I^0(\rest{aL+(b-c)A}{bY})$ 
induced by the surjective homomorphism
$H^0(aL+(b-c)A) \to I^0(\rest{aL+(b-c)A}{bY})$ 
and the $L^2$-norm $\Vert\cdot\Vert_{L^2}^{a\overline{L}+(b-c)\overline{A}}$
of $H^0(aL + (b-c)A)$.
Note that $I^0(\rest{aL+(b-c)A}{bY})$ is torsion free because
$bY$ is flat over $\ZZ$.

\begin{Claim}
\label{claim:thm:h:0:estimate:big:3}
For all integers $a, b, c$ with $a \geq b \geq c \geq 0$ and $a \geq 2$,
\begin{multline*}
\ah\left(H^0(aL+ (b-c)A), \Vert\cdot\Vert_{L^2}^{a\overline{L}+ (b-c)\overline{A}}\right) \leq
\ah\left(H^0(aL - cA), \Vert\cdot\Vert_{L^2, s^b, \sub}^{a\overline{L}+(b-c)\overline{A}}\right)   \\
+ \ah\left(I^0(\rest{aL+(b-c)A}{bY}), \Vert\cdot\Vert_{L^2,\quot}^{a\overline{L}+(b-c)\overline{A}}\right) 
+ C_2 a^{d-1}\log(a).
\end{multline*}
\end{Claim}

Using an exact sequence
\[
0 \to H^0(aL - cA) \overset{s^b}{\longrightarrow} H^0(aL+(b-c)A) \to I^0(\rest{aL+(b-c)A}{bY}) \to 0,
\]
we have a normed exact sequence
\begin{multline*}
0 \to \left(H^0(aL-cA), \Vert\cdot\Vert_{L^2, s^b, \sub}^{a\overline{L}+(b-c)\overline{A}}\right) 
\to
\left(H^0(a L+(b-c)A), \Vert\cdot\Vert_{L^2}^{a\overline{L}+(b-c)\overline{A}}\right) \\
\to \left(I^0(\rest{aL+(b-c)A}{mY}), \Vert\cdot\Vert_{L^2,\quot}^{a\overline{L}+(b-c)\overline{A}}\right) \to 0,
\end{multline*}
where $ \Vert\cdot\Vert_{L^2, s^b, \sub}^{a\overline{L}+(b-c)\overline{A}}$ 
is the subnorm of $H^0(aL-cA)$ induced by the injective homomorphism
$H^0(aL - cA) \overset{s^b}{\longrightarrow} H^0(aL+(b-c)A)$ and
the $L^2$-norm $\Vert\cdot\Vert_{L^2}^{a\overline{L}+(b-c)\overline{A}}$
of $H^0(a L+(b-c)A)$.
Thus, by \eqref{eqn:rem:prop:basic:property:h:0:h:1:x:4}, it yields the claim because
\[
\rank H^0(a L-cA) \leq \rank H^0(a(L+A)).
\]

\medskip
Next we claim the following:

\begin{Claim}
\label{claim:thm:h:0:estimate:big:4}
There are constants $a_0$ and $C_5$ depending only on
$\overline{L}$ and $\overline{A}$ such that
\begin{multline*}
\ah\left(H^0(aL - cA), \Vert\cdot\Vert_{L^2, s^b,\sub}^{a\overline{L}+(b-c)\overline{A}}\right) \leq
\ah\left(H^0(aL - cA), \Vert\cdot\Vert_{L^2}^{a\overline{L} - c \overline{A}}\right) \\
+ C_5 b a^{d-1} + C_2 a^{d-1}\log(a).
\end{multline*}
for all integers $a, b, c$ with $a \geq b \geq c \geq 0$ and $a \geq a_0$.
\end{Claim}

Note that $\Vert\cdot\Vert_{L^2, s^b,\sub}^{a\overline{L}+(b-c)\overline{A}} \leq
 \Vert\cdot\Vert_{L^2}^{a\overline{L} - c \overline{A}}$. Thus,
 by \eqref{eqn:rem:prop:basic:property:h:0:h:1:x:2},
\begin{multline*}
\ah\left(H^0(aL-cA), \Vert\cdot\Vert_{L^2}^{a\overline{L}- c\overline{A}}\right) - 
\ah\left(H^0(aL-cA), \Vert\cdot\Vert_{L^2, s^b,\sub}^{a\overline{L}+(b-c)\overline{A}}\right) + C_2 a^{d-1} \log(a) \\
\geq \achi\left(H^0(aL-cA), \Vert\cdot\Vert_{L^2}^{a\overline{L} - c\overline{A}}\right) - 
\achi\left(H^0(aL-cA), \Vert\cdot\Vert_{L^2, s^b,\sub}^{a\overline{L}+(b-c)\overline{A}}\right).
\end{multline*}
Therefore it is sufficient to find positive constants $a_0$ and $C_5$ such that
\[
\achi\left(H^0(aL-cA), \Vert\cdot\Vert_{L^2}^{a\overline{L}- c\overline{A}}\right) - 
\achi\left(H^0(aL-cA), \Vert\cdot\Vert_{L^2, s^b, \sub}^{a\overline{L}+(b-c)\overline{A}}\right)
\geq -C_5 b a^{d-1}
\]
for all $a, b, c$ with $a \geq b \geq c \geq 0$ and $a \geq a_0$.
This is nothing more than a consequence of Corollary~\ref{cor:comp:sub:L2:ball}.

\medskip
Let $k$ be an integer with $0 \leq k < b$.
Let $\Vert\cdot\Vert_{L^2, s^k, \sub,\quot}^{a\overline{L}+(b-c)\overline{A}}$ be
the quotient norm of $I^0(Y, \rest{aL + (b-c-k)A}{Y})$ induced by a surjective homomorphism
\[
H^0(aL + (b-c-k)A) \to I^0(Y, \rest{aL + (b-c-k)A}{Y})
\]
and
$\Vert\cdot\Vert_{L^2, s^k, \sub}^{a\overline{L}+(b-c)\overline{A}}$ of
$H^0(aL + (b-c-k)A)$.

\begin{Claim}
\label{claim:thm:h:0:estimate:big:5}
There is a constant $C_6$ and $C_7$ depending only on $\overline{L}$ and
$\overline{A}$ such that
\[
\ah\left(I^0(Y, \rest{aL + (b-c-k)A}{Y}), 
\Vert\cdot\Vert_{L^2, s^k, \sub,\quot}^{a\overline{L}+(b-c)\overline{A}}\right)
\leq C_6 a^{d-1} + C_7 a^{d-2}\log(a)
\]
for all integers $a, b, c, k$ with $a \geq b \geq c \geq 0$, $a \geq 2$ and $0 \leq k < b$.
\end{Claim}

Let us choose a small open set $U$ of $X(\CC)$ such that
the closure of $U$ does not meet with $Y(\CC)$
and $U$ is not empty on each connected component of $X(\CC)$.
Then, applying Lemma~\ref{lem:comp:norm:X:U} to the cases
$L_{\CC}, A_{\CC}$ and $L_{\CC}, -A_{\CC}$,
there are  constant $D_1 \geq 1$ and $D'_1 \geq 1$ such that
\[
 D'_1 D_1^{l + \vert m \vert} 
 \int_U \vert u \vert^2 \Omega \geq \int_{X(\CC)} \vert u \vert^2 \Omega
\]
for all integers $l, m$ with $l  \geq 0$ and all $u \in H^0(X(\CC), lL + mA)$.
Since $0 < \inf_{x \in U} \{ \vert s \vert(x)\} <1$,
if we set
\[
 D_2 = 1/ \inf_{x \in U} \{ \vert s \vert(x)\},
\]
then $D_2 > 1$.
Thus, if we set $D_3 = \max\{ D_2, D_1\}$, then,
for $u \in H^0(X, aL +(b-c-k)A)$,
\begin{align*}
\int_{X(\CC)}  \vert s^k \otimes u \vert^2 \Omega & 
\geq \int_U \vert s^k \otimes u \vert ^2 \Omega \geq D_2^{-2k} \int_U \vert u \vert^2 \Omega \\
& \geq  D_2^{-2k} {D'_1}^{-1} D_1^{-(a+\vert b-c-k\vert )} \int_{X(\CC)} \vert u \vert^2 \Omega
\geq {D'_1}^{-1} D_3^{-4a}  \int_{X(\CC)} \vert u \vert^2 \Omega,
\end{align*}
which means that
\[
\Vert\cdot\Vert_{L^2, s^k,\sub}^{a\overline{L}+(b-c)\overline{A}} \geq {D'_1}^{-1/2} D_3^{-2a}
\Vert\cdot\Vert_{L^2}^{a\overline{L}+(b-c-k)\overline{A}}.
\]
Hence
\[
\Vert\cdot\Vert_{L^2, s^k,\sub,\quot}^{a\overline{L}+(b-c)\overline{A}} \geq {D'_1}^{-1/2} D_3^{-2a}
\Vert\cdot\Vert_{L^2,\quot}^{a\overline{L}+(b-c-k)\overline{A}},
\]
where $\Vert\cdot\Vert_{L^2,\quot}^{a\overline{L}+(b-c-k)\overline{A}}$ is the quotient norm of
$I^0(Y, \rest{aL + (b - c -k)A}{Y})$ induced by a surjective homomorphism
\[
H^0(X, aL + (b -c-k)A) \to I^0(Y, \rest{aL + (b - c -k)A}{Y}).
\]
Note that $e^x \geq x + 1$ for $x \geq 0$.
Thus, applying Corollary~\ref{cor:norm:comp:L:2:subvariety} to the cases
$L_{\CC}, A_{\CC}$ and $L_{\CC}, -A_{\CC}$,
there are constants $D_4, D'_4  \geq 1$ such that
\begin{multline*}
 \Vert\cdot\Vert_{L^2,\quot}^{a\overline{L}+(b-c-k)\overline{A}} \geq
{D'_4}^{-1/2} D_4^{-(a+\vert b-c - k\vert)/2}\Vert\cdot\Vert_{L^2}^{\rest{a\overline{L}+(b-c-k)\overline{A}}{Y}} \\
\geq 
{D'_4}^{-1/2} D_4^{-a}\Vert\cdot\Vert_{L^2}^{\rest{a\overline{L}+(b-c-k)\overline{A}}{Y}}
\end{multline*}
on $I^0(Y, \rest{aL+(b-c-k)A}{Y})$,
where the volume form on $Y$ is
given by the $C^{\infty}$-hermitian invertible sheaf
$\rest{\overline{L} + \overline{A}}{Y}$.
Therefore, if we set $D_5 = \max \{ D_3, D_4 \}$ and $D'_5 = \max \{ D'_1, D'_4 \}$, then
\[
\Vert\cdot\Vert_{L^2, s^k,\sub,\quot}^{a\overline{L}+(b-c)\overline{A}} \geq {D'_5}^{-1}D_5^{-3a} 
\Vert\cdot\Vert_{L^2}^{\rest{a\overline{L}+(b-c-k)\overline{A}}{Y}}
\]
on $I^0(Y, \rest{aL + (b-c-k)A}{Y})$. Thus, 
by \eqref{eqn:rem:prop:basic:property:h:0:h:1:x:3},
\begin{multline*}
\ah\left(I^0(Y, \rest{aL + (b-c-k)A}{Y}), 
\Vert\cdot\Vert_{L^2, s^k,\sub,\quot}^{a\overline{L}+(b-c)\overline{A}}\right) \\
\leq \ah\left(I^0(Y, \rest{aL + (b-c-k)A}{Y}),
\Vert\cdot\Vert_{L^2}^{\rest{a\overline{L}+(b-c-k)\overline{A}}{Y}}\right) \\
\qquad\qquad\qquad\qquad+ \log({D'_5}D_5^{3a} )C_3 a ^{d-2}  + C_4 a^{d-2}\log(a) \\
\leq \ah\left(H^0(Y, \rest{aL+(b-c-k)A}{Y}),
\Vert\cdot\Vert_{L^2}^{\rest{a\overline{L}+(b-c-k)\overline{A}}{Y}}\right)\\
\qquad\qquad\qquad\qquad+\log({D'_5}D_5^{3a} )C_3 a ^{d-2}  + C_4 a^{d-2}\log(a).
\end{multline*}
Let $\widetilde{Y}$ be the normalization of $Y$. 
Let $t$ be a non-zero section as in Claim~\ref{claim:thm:h:0:estimate:big:2}.
Then  $t$ gives rise to a relation $\rest{\overline{L}}{\widetilde{Y}} \leq \rest{n\overline{A}}{\widetilde{Y}}$
(cf. Conventions and terminology~\ref{CT:Arakelov:order:hermitian}).
Thus 
\[
\rest{a\overline{L}+(b-c-k)\overline{A}}{\widetilde{Y}} \leq \rest{(an + b-c-k)\overline{A}}{\widetilde{Y}}.
\]
Therefore,
\begin{multline*}
 \ah\left(H^0(Y, \rest{aL+(b-c-k)A}{Y}),
\Vert\cdot\Vert_{L^2}^{\rest{a\overline{L}+(b-c-k)\overline{A}}{Y}}\right)\\
\leq
\ah\left(H^0(\widetilde{Y}, \rest{aL+(b-c-k)A}{\widetilde{Y}}),
\Vert\cdot\Vert_{L^2}^{\rest{a\overline{L}+(b-c-k)\overline{A}}{\widetilde{Y}}} \right) \\
\leq 
\ah\left(H^0(\widetilde{Y}, \rest{(an + b-c-k)A}{\widetilde{Y}}),
\Vert\cdot\Vert_{L^2}^{\rest{(an + b-c-k)\overline{A}}{\widetilde{Y}}}\right).
\end{multline*}
Further,
by Lemma~\ref{lem:upper:estimate:h:0},
there is a positive constant $D_6$
with
\[
\ah\left(H^0(\widetilde{Y}, \rest{nA}{\widetilde{Y}}),
\Vert\cdot\Vert_{L^2}^{n \rest{\overline{A}}{\widetilde{Y}}}\right)
\leq D_6 n^{d-1}
\]
for all $n \geq 1$. Thus the claim follows.

\medskip
Finally we claim the following:

\begin{Claim}
\label{claim:thm:h:0:estimate:big:6}
There is a constant $C_7$ depending only on $\overline{L}$ and
$\overline{A}$ such that
\[
\ah\left(I^0(\rest{(aL + (b-c)A)}{bY}), \Vert\cdot\Vert_{L^2,\quot}^{a\overline{L} + (b-c)\overline{A}}\right) \\
\leq C_6 b a^{d-1} + (C_4 + C_7) a^{d-1}\log(a)
\]
for all integers $a, b, c$ with $a \geq b \geq c \geq 0$, $a \geq 2$.
\end{Claim}

A commutative diagram
\[
\begin{CD}
0 @>>> -(k+1)A @>{s^{k+1}}>> \OO_X @>>> \OO_{(k+1)Y} @>>> 0 \\
@. @VV{s}V @| @VVV @. \\
0 @>>> -kA @>{s^k}>> \OO_X @>>> \OO_{kY} @>>> 0 \\
@. @VVV @. @. @. \\
  @. \rest{-kA}{Y} @.   @.  @.
\end{CD}
\]
yields an injective homomorphism $\alpha_k : \rest{-kA}{Y} \to \OO_{(k+1)Y}$
together with a commutative diagram
\[
\begin{CD}
0 @>>> -kA @>{s^k}>> \OO_X @>>> \OO_{kY} @>>> 0 \\
@. @VVV @VVV @| @. \\
0 @>>> \rest{-kA}{Y} @>{\alpha_k}>> \OO_{(k+1)Y} @>>> \OO_{kY} @>>> 0,
\end{CD}
\]
where two horizontal sequences are exact.
Thus, tensoring the above diagram with $aL + (b-c)A$,
we have the following commutative diagram:
\[
{
\xymatrix{
0 \ar[r] & aL + (b-c-k)A \ar[r]^{s^k} \ar[d] & aL + (b-c)A \ar[r] \ar[d] & \rest{aL + (b-c)A}{kY} \ar[r] \ar@{=}[d] & 0 \\
0 \ar[r] & \rest{aL + (b-c-k)A}{Y} \ar[r]^{\alpha_k} & \rest{aL + (b-c)A}{(k+1)Y} \ar[r] & \rest{aL + (b-c)A}{kY} \ar[r] & 0
}}
\]
Therefore we have an exact sequence
\begin{multline*}
0 \to I^0(\rest{(aL + (b-c-k)A)}{Y}) \to
I^0(\rest{(aL + (b-c)A)}{(k+1)Y}) \\
\to I^0(\rest{(aL + (b-c)A)}{kY}) \to 0
\end{multline*}
Note that in the commutative diagram
\[
\begin{CD}
H^0(aL + (b-c-k)A) @>{s^k}>> H^0(aL + (b-c)A) \\
@VVV @VVV \\
I^0(\rest{(aL + (b-c-k)A)}{Y}) @>{\alpha_k}>>  I^0( \rest{(aL + (b-c)A)}{(k+1)Y}),
\end{CD}
\]
the two vertical arrows have the same kernel.
Thus, by Lemma~\ref{lem:sub:quot:vs:quot:sub},
\begin{multline*}
0 \to
\left(I^0(\rest{(aL + (b-c-k)A)}{Y}), 
\Vert\cdot\Vert_{L^2, s^k, \sub,\quot}^{a\overline{L} + (b-c)\overline{A}}\right) \\
\to \left(I^0(\rest{(aL + (b-c)A)}{(k+1)Y}), 
\Vert\cdot\Vert_{L^2,\quot}^{a\overline{L} + (b-c)\overline{A}}\right) \\
\to \left(I^0(\rest{(aL + (b-c)A)}{kY}), 
\Vert\cdot\Vert_{L^2,\quot}^{a\overline{L} + (b-c)\overline{A}}\right) \to 0
\end{multline*}
is a normed exact sequence, where
for each $1 \leq i \leq b$, the norm
$\Vert\cdot\Vert_{L^2,\quot}^{a\overline{L} + (b-c)\overline{A}}$ of
$I^0(\rest{(aL + (b-c)A)}{iY})$
is the quotient norm
induced by the surjective homomorphism
$H^0(aL + (b-c)A) \to I^0(\rest{(aL + (b-c)A)}{iY})$ and the $L^2$-norm
$\Vert\cdot\Vert_{L^2}^{a\overline{L} + (b-c)\overline{A}}$ of
$H^0(aL + (b-c)A)$.
Therefore, 
by \eqref{eqn:rem:prop:basic:property:h:0:h:1:x:4},
\begin{multline*}
\ah\left(I^0(\rest{(aL + (b-c)A)}{(k+1)Y}), 
\Vert\cdot\Vert_{L^2,\quot}^{a\overline{L} + (b-c)\overline{A}}\right) \\
- \ah\left(I^0(\rest{(aL + (b-c)A)}{kY}), 
\Vert\cdot\Vert_{L^2,\quot}^{a\overline{L} + (b-c)\overline{A}}\right) \\
\leq \ah\left(I^0(\rest{(aL + (b-c-k)A)}{Y}), 
\Vert\cdot\Vert_{L^2, s^k, \sub,\quot}^{a\overline{L} + (b-c)\overline{A}}\right) +C_4 a^{d-2}\log(a).
\end{multline*}
Thus, taking $\sum_{k=1}^{b-1}$, the above yields
\begin{multline*}
\ah\left(I^0(\rest{(aL + (b-c)A)}{bY}), 
\Vert\cdot\Vert_{L^2,\quot}^{a\overline{L} + (b-c)\overline{A}}\right) \\
\leq \sum_{k=0}^{b-1}
\ah\left(I^0(\rest{(aL + (b-c-k)A)}{Y}), 
\Vert\cdot\Vert_{L^2, s^k, \sub,\quot}^{a\overline{L} + (b-c)\overline{A}}\right) \\
+ (b-1)C_4 a^{d-2}\log(a).
\end{multline*}
Therefore, using Claim~\ref{claim:thm:h:0:estimate:big:5},
we have the claim.

\medskip
Gathering Claim~\ref{claim:thm:h:0:estimate:big:3}, Claim~\ref{claim:thm:h:0:estimate:big:4}
and Claim~\ref{claim:thm:h:0:estimate:big:6},
if we set $C = C_5 + C_6$ and $D = 2C_2 + C_4 + C_7$, then
\begin{multline*}
\ah\left(H^0(aL + (b-c)A),
\Vert\cdot\Vert^{a\overline{L} + (b-c)\overline{A}}_{L^2}\right) \\
\leq \ah\left(H^0(a L - c A), 
\Vert\cdot\Vert^{a \overline{L} - c \overline{A}}_{L^2}\right)
+ C b a^{d-1} + D a^{d-1} \log(a)
\end{multline*}
for all $a  \geq b \geq c \geq 0$ and $a \geq a_0$.

\medskip
Finally let us consider the sup-version of our estimate.
First of all, since 
\[
\Vert\cdot\Vert^{a\overline{L} + (b-c)\overline{A}}_{\sup}
\geq \Vert\cdot\Vert^{a\overline{L} + (b-c)\overline{A}}_{L^2},
\]
we have
\[
\ah\left(H^0(aL + (b-c)A),
\Vert\cdot\Vert^{a\overline{L} + (b-c)\overline{A}}_{\sup}\right) \leq
\ah\left(H^0(aL + (b-c)A),
\Vert\cdot\Vert^{a\overline{L} + (b-c)\overline{A}}_{L^2}\right).
\]
Moreover, by virtue of Gromov's inequality, there is a constant $C_8 \geq 1$
\[
\Vert\cdot\Vert^{a \overline{L} - c \overline{A}}_{L^2} \geq
C_8^{-1} (a+c+1)^{-(d-1)} \Vert\cdot\Vert^{a \overline{L} - c \overline{A}}_{\sup}
\]
for all $a, c \geq 0$.
Thus, since $a \geq c$,
\begin{multline*}
\ah\left(H^0(a L - c A), 
\Vert\cdot\Vert^{a \overline{L} - c \overline{A}}_{L^2}\right) \\
\leq \ah\left(H^0(a L - c A), 
C_8^{-1} (a+c+1)^{-(d-1)} \Vert\cdot\Vert^{a \overline{L} - c \overline{A}}_{\sup}\right) \\
\leq \ah\left(H^0(a L - c A), 
\Vert\cdot\Vert^{a \overline{L} - c \overline{A}}_{\sup}\right) \\
+  \log(C_8(2a + 1)^{d-1}) C_1 a^{d-1} + C_2 a^{d-1} \log(a)
\end{multline*}
for all $a \geq c \geq 0$.
Therefore we obtain the sup-version.
\QED

Let $R$ be an integral domain such that
$R$ is flat and finite over $\ZZ$.
Let $K$ be a quotient field of $R$.
Note that $K$ is a number field.
Let $K(\CC)$ be the set of all embeddings $K \hookrightarrow \CC$ of fields.
Let $L$ be a finitely generated and free $R$-module of rank $1$.
For each $\sigma \in K(\CC)$,
the tensor product $L \otimes_{R} \CC$ in terms of
the embedding $\sigma : K \hookrightarrow \CC$ is
denoted by $L_{\sigma}$. For each $\sigma \in K(\CC)$,
let $\vert\cdot\vert_{\sigma}$ be a norm of $L_{\sigma}$.
The collection $\left(L, \{ \vert\cdot\vert_{\sigma} \}_{\sigma \in K(\CC)}\right)$
is called a normed invertible $R$-module.
For simplicity,
$\left(L, \{ \vert\cdot\vert_{\sigma} \}_{\sigma \in K(\CC)}\right)$ is often
denoted by $(L, \vert\cdot\vert)$ or $\overline{L}$.
We define $\Vert\cdot\Vert_{\sup}^{\overline{L}}$ by
\[
\Vert s \Vert_{\sup}^{\overline{L}} = \max \{ \vert s \vert_{\sigma} \mid \sigma \in K(\CC) \}.
\]
Then $(L, \Vert\cdot\Vert^{\overline{L}}_{\sup})$ is a normed
finitely generated free $\ZZ$-module.

\begin{Proposition}
\label{prop:0:estimate:big:curve}
Let $\overline{L}$ and $\overline{A}$ be normed invertible $R$-modules of rank $1$.
We assume that there is $s \in A$ with $s \not= 0$ and $\Vert s \Vert_{\sup}^{\overline{A}} \leq 1$.
Then there are positive constants $C$ and $D$ depending only on
$\overline{L}$ and $\overline{A}$ such that
\begin{multline*}
\ah\left(a L + (b-c)A,
\Vert\cdot\Vert^{a \overline{L}+(b-c)\overline{A}}_{\sup}\right)
\leq
\ah\left(a L - c A, 
\Vert\cdot\Vert^{a \overline{L} - c \overline{A}}_{\sup}\right)
+ C b + D
\end{multline*}
for all non-negative integers $a, b, c$.
\end{Proposition}

\Proof
Let $\Vert\cdot\Vert_{\sup, s^b, \sub}^{a \overline{L} + (b-c) \overline{A}}$
be the subnorm of $a L - c A$ induced by the injective homomorphism
$a L - c A \overset{s^{b}}{\longrightarrow} a L +(b-c) A$
and the norm $\Vert\cdot\Vert_{\sup}^{a \overline{L} + (b-c) \overline{A}}$
of $a L +(b-c) A$. 
Then, by (4) of Proposition~\ref{prop:basic:property:h:0:h:1:x},
we have
\begin{multline*}
\ah\left(a L +(b-c)A, 
\Vert\cdot\Vert_{\sup}^{a \overline{L} + (b-c) \overline{A}} \right)
\leq \ah\left( a L - c A,
\Vert\cdot\Vert_{\sup, s^b, \sub}^{a \overline{L} + (b-c) \overline{A}} \right) \\
+ \log \#(\Coker(a L-c A \to
a L + (b-c) A)) + D,
\end{multline*}
where $d = [K : \QQ]$ and $D = (\log (18) + 2) (d + 1)\log(d + 1)$.
Note that
\begin{multline*}
\log  \#(\Coker(a L - c A \overset{s^{b}}{\longrightarrow}
a L +(b-c) A)) \\
=
\log  \#(\Coker(R \overset{s^{b}}{\longrightarrow} bA) \otimes (a L - c A)) \\
= \log  \#(\Coker(R \overset{s^{b}}{\longrightarrow} bA)).
\end{multline*}
Let us consider a sequence of injective homomorphisms:
\[
R \overset{s}{\longrightarrow}
A \overset{s}{\longrightarrow} \cdots
 \overset{s}{\longrightarrow} b A.
\]
Then
\begin{multline*}
\log  \#(\Coker(R \overset{s^{b}}{\longrightarrow}
bA)) = \sum_{i=1}^b
\log \#(\Coker((i-1)A \overset{s}{\longrightarrow} iA)) \\
=
b \cdot \log \#(\Coker(R \overset{s}{\longrightarrow} A)).
\end{multline*}
On the other hand,
for all $t \in a L - c A$,
\[
\Vert s^b \otimes t \Vert_{\sup}^{a\overline{L}+(b-c)\overline{A}}
\geq \left( \min \{ \vert s \vert_{\sigma} \mid \sigma \in K(\CC) \} \right)^b
\Vert t \Vert_{\sup}^{a \overline{L} - c \overline{A}}.
\]
Thus, by (3) of Proposition~\ref{prop:basic:property:h:0:h:1:x},
\begin{multline*}
 \ah\left( a L - c A,
\Vert\cdot\Vert_{\sup, s^b,\sub}^{a \overline{L}+(b-c) \overline{A}} \right)
\leq 
 \ah\left( a L - c A,
\Vert\cdot\Vert_{\sup}^{a \overline{L} - c \overline{A}} \right) 
+  b \log(C'(s)) d+ D,
\end{multline*}
where $C'(s) = \min \{ \vert s \vert_{\sigma} \mid \sigma \in K(\CC) \}$.
Therefore,
\begin{multline*}
\ah\left(a L + (b-c) A, 
\Vert\cdot\Vert_{\sup}^{a \overline{L} +(b-c) \overline{A}} \right)
\leq \ah\left( a L - c A,
\Vert\cdot\Vert_{\sup}^{a \overline{L} - c \overline{A}} \right) \\
+ b( \log \#(\Coker(R \overset{s}{\longrightarrow} A)) +  \log(C'(s)))d + 2D.
\end{multline*}
\QED

Finally we consider the following  lemma which guarantees the existence of a good
$C^{\infty}$-hermitian invertible sheaf $\overline{A}$ satisfying the assumptions (i), (ii) and (iii) of
Theorem~\ref{thm:h:0:estimate:big}.

\begin{Lemma}
\label{lem:ample:generic:smooth}
Let $X$ be a projective and generically smooth arithmetic variety
of dimension $d \geq 2$, and let
$\overline{A}$ be an ample  $C^{\infty}$-hermitian invertible sheaf on $X$.
Then, for any $C^{\infty}$-hermitian invertible sheaf $\overline{L}$ on $X$,
there is a positive integer $n_0$ such that,
for all $n \geq n_0$,
$n \overline{A}$ satisfies the assumptions \rom{(i)}, \rom{(ii)} and \rom{(iii)} of
Theorem~\rom{\ref{thm:h:0:estimate:big}}.
\end{Lemma}

\Proof
This is a consequence of arithmetic Bertini's theorem (cf. \cite{MoBogo2}).
We can give however an easy and direct proof of the lemma as follows:
It is easy to find $n_0$ for the assumptions (i) and (ii).
In addition to (i) and (ii),
we choose $n_0$ such that $nA$ is very ample for all $n \geq n_0$.
Let $\pi : X \to \Spec(\ZZ)$ be the structure morphism and
$S$ the minimal finite set of $\Spec(\ZZ) \setminus \{ 0 \}$ such that
$\pi^{-1}(\Spec(\ZZ) \setminus S)$ is regular.
Let $Z_1, \ldots, Z_r$ be all irreducible components of
$\pi^{-1}(S)$, and let $x_1, \ldots, x_r$ be closed points of $X$
with $x_i \in Z_i$ for all $i$.
Let $m_1, \ldots, m_r$ be the maximal ideals corresponding to
$x_1, \ldots, x_r$.
Then there is a positive integer $n_1$ such that,
for all $n \geq n_1$,
$H^1(X, nA \otimes m_1 \cdots m_r) = 0$, which means that
the natural homomorphism
\[
H^0(X, nA) \to \bigoplus_{i=1}^n nA \otimes (\OO_X/m_i)
\]
is surjective. 
Thus if $n \geq \max \{ n_0, n_1 \}$,
then $nA$ is very ample and there is
a non-zero section $t_n$ of $H^0(X, nA)$
with $t_n(x_i) \not= 0$ for all $x_i$.
We set $\gamma(s) = t _n+ l s$ for $s \in H^0(X, nA)$,
where $l = \prod_{s \in S} \characteristic(\kappa(s))$
and $\kappa(s)$ is the residue field of $\ZZ$ at $s$.
Note that $\gamma(s)(x_i) \not= 0$ for all $i$.
In particular, every vertical component of $\zeros(\gamma(s))$
is contained $\pi^{-1}(\Spec(\ZZ) \setminus S)$.
On the other hand, it is easy to see that the set
$\{ \gamma(s) \mid s \in H^0(X, nA)\}$ is Zariski dense
in a vector space $H^0(X_{\QQ}, nA_{\QQ}) = H^0(X, nA) \otimes \QQ$.
Thus, by Bertini's theorem,
there is $s \in H^0(X, nA)$ such that
$\zeros(\gamma(s))$ is smooth over $\QQ$.
\QED

\section{Volume function for $C^{\infty}$-hermitian invertible sheaves \\ and its basic properties}

Let $X$ be a projective arithmetic variety
of dimension $d$.
For a $C^{\infty}$-hermitian invertible sheaf $\overline{L}$ on $X$,
the {\em arithmetic volume of $\overline{L}$} is
defined by
\[
\avol(\overline{L}) = \limsup_{m \to \infty}
\frac{\ah(H^0(X, mL), \Vert\cdot\Vert_{\sup}^{m \overline{L}})}{m^d/d!}.
\]
This number is a finite real number by Lemma~\ref{lem:upper:estimate:h:0}.
Moreover, if $\overline{L}$ is ample, then
\[
\avol(\overline{L}) = \adeg(\acherncl_1(\overline{L})^{\cdot d}).
\]

First let us consider elementary properties of volume function:

\begin{Proposition}
\label{prop:basic:properties:volume}
Let $\overline{L}$ and $\overline{M}$ be $C^{\infty}$-hermitian invertible sheaves on $X$.
Then we have the following:
\begin{enumerate}
\renewcommand{\labelenumi}{(\arabic{enumi})}
\item
If $\overline{L} \leq \overline{M}$ 
\rom{(}Conventions and terminology~\rom{\ref{CT:Arakelov:order:hermitian}}\rom{)},
then $\avol(\overline{L}) \leq \avol(\overline{M})$.

\item
Let $\vert\cdot\vert_L$ be the hermitian norm of $\overline{L}$.
For a real number $\lambda$, we set 
\[
\overline{L}^{\lambda} = (L, \exp(-\lambda)\vert\cdot\vert_L).
\]
If $\lambda \geq 0$, then we have
\[
\begin{cases}
\avol(\overline{L}) \leq
\avol(\overline{L}^{\lambda}) \leq \avol(\overline{L}) +d \lambda\vol(L_{\QQ}), \\
\avol(\overline{L}) - d \lambda \vol(L_{\QQ}) \leq
\avol(\overline{L}^{- \lambda}) \leq \avol(\overline{L}),
\end{cases}
\]
where $\vol(L_{\QQ})$ is the geometric volume of $L_{\QQ}$ on $X_{\QQ}$.

\item
${\displaystyle
\avol(\overline{L}) = \limsup_{m \to \infty} \frac{\log \# \{ s \in H^0(X, mL) \mid \Vert s \Vert^{m\overline{L}}_{\sup} < 1 \}}{m^d/d!}.
}$
%
\end{enumerate}
\end{Proposition}

\Proof
(1) Since $\overline{L} \leq \overline{M}$, we have $m \overline{L} \leq m \overline{M}$ for
all $m \geq 1$.
Thus 
\[
\ah\left(H^0(X, mL), \Vert\cdot\Vert_{\sup}^{m \overline{L}} \right) \leq
 \ah\left(H^0(X, mM), \Vert\cdot\Vert_{\sup}^{m \overline{M}}\right)
\]
for all $m \geq 1$. Hence $\avol(\overline{L}) \leq \avol(\overline{M})$.

\medskip
(2)
Since $\Vert\cdot\Vert^{m\overline{L}^{\lambda}}_{\sup} = \exp(-m\lambda)
\Vert\cdot\Vert^{m\overline{L}}_{\sup}$, by using  \eqref{eqn:rem:prop:basic:property:h:0:h:1:x:3},
there is a positive constant $C$ such that 
\begin{multline*}
0 \leq \ah(H^0(X, mL), \Vert\cdot\Vert^{\overline{L}^{\lambda}}_{\sup}) -
\ah(H^0(X, mL), \Vert\cdot\Vert^{\overline{L}}_{\sup}) \\
\leq \lambda m \dim H^0(X_{\QQ}, mL_{\QQ}) + C m^{d-1}\log(m)
\end{multline*}
for $m \gg 1$. Thus we obtain the first inequalities.
These implies that
\[
\avol(\overline{L}^{-\lambda}) \leq
\avol\left(\left(\overline{L}^{-\lambda}\right)^{\lambda}\right) \leq \avol(\overline{L}^{-\lambda}) +d \lambda\vol(L_{\QQ}),
\]
which is nothing more than the second inequalities
because
$\left(\overline{L}^{-\lambda}\right)^{\lambda} = \overline{L}$.

\medskip
(3)
For a positive real number $\lambda$,
\begin{multline*}
\aH\left(H^0(X, mL), \Vert\cdot\Vert_{\sup}^{m\overline{L}^{-\lambda}}\right) \\
\subseteq
\left\{ s \in H^0(X, mL) \mid \Vert s \Vert^{m\overline{L}}_{\sup} < 1 \right\} \\
\subseteq
\aH\left(H^0(X, mL), \Vert\cdot\Vert_{\sup}^{m\overline{L}}\right)
\end{multline*}
because $\Vert\cdot\Vert_{\sup}^{m\overline{L}^{-\lambda}} = \exp(m\lambda) \Vert\cdot\Vert_{\sup}^{m\overline{L}}$.
Thus, using (2), we have
\[
\avol(\overline{L}) - d \lambda \vol(L_{\QQ}) \leq
 \limsup_{m \to \infty} \frac{\log \# \{ s \in H^0(X, mL) \mid \Vert s \Vert^{m\overline{L}}_{\sup} < 1 \}}{m^d/d!}
 \leq \avol(\overline{L}),
\]
which shows the assertion
because $\lambda$ is an arbitrary positive number.
%
\QED

The following theorem shows that the volume function is a birational invariant.

\begin{Theorem}
\label{thm:invariance:vol:birat:morphism}
Let $\pi : X' \to X$ be a birational morphism of projective arithmetic varieties, and
let  $\overline{L}$ and $\overline{N}$ be $C^{\infty}$-hermitian invertible sheaves on $X$.
Then
\begin{multline*}
\limsup_{m \to \infty} 
\frac{\ah\left(H^0(X, mL + N), \Vert\cdot\Vert_{\sup}^{m\overline{L} + \overline{N}}\right)}{m^d} \\
=
\limsup_{m \to \infty} 
\frac{\ah\left(H^0(X', \pi^*(mL + N)), \Vert\cdot\Vert_{\sup}^{\pi^*(m\overline{L} + \overline{N})}\right)}{m^d}.
\end{multline*}
In particular,
$\avol(\overline{L}) = \avol(\pi^*(\overline{L}))$.
\end{Theorem}

\Proof
The proof of this theorem is similar to \cite[Theorem~2.2]{Yuan}.
First of all, note that
\begin{multline*}
\limsup_{m \to \infty} 
\frac{\ah\left(H^0(X, mL + N), \Vert\cdot\Vert_{\sup}^{m\overline{L} + \overline{N}}\right)}{m^d} \\
\leq
\limsup_{m \to \infty} 
\frac{\ah\left(H^0(X', \pi^*(mL + N)), \Vert\cdot\Vert_{\sup}^{\pi^*(m\overline{L} + \overline{N})}\right)}{m^d}.
\end{multline*}
Thus, considering a generic resolution of singularities of $X'$, we may assume that
$X'$ is generically smooth.

Let us consider an exact sequence:
\[
0 \to mL + N \to \pi_*(\pi^*(mL + N)) \to (mL + N)\otimes (\pi_*(\OO_{X'})/\OO_X) \to 0.
\]
The image of the natural homomorphism
\[
H^0(X', \pi^*(mL+ N)) \to H^0(X,  (mL + N)\otimes (\pi_*(\OO_{X'})/\OO_X))
\]
is denoted by $\Gamma(X'/X, mL+ N)$.
Let $\Vert\cdot\Vert^{\pi^*(m\overline{L}+\overline{N})}_{\sup, \quot}$ be
the quotient norm of $\Gamma(X'/X, mL+ N)$
induced by the surjective homomorphism 
\[
H^0(X', \pi^*(mL+ N)) \to \Gamma(X'/X, mL+N)
\]
and
the sup-norm $\Vert\cdot\Vert^{\pi^*(m\overline{L}+\overline{N})}_{\sup}$
of $H^0(X', \pi^*(mL+N))$.

\begin{Claim}
\label{claim:prop:invariance:vol:birat:morphism:1}
If $\pi$ is finite and $\overline{L}$ is ample, then
\[
\ah\left( \Gamma(X'/X, mL+N), \Vert\cdot\Vert^{\pi^*(m\overline{L}+\overline{N})}_{\sup,\quot} \right)
\leq o(m^d).
\]
\end{Claim}

We fix a normalized volume form $\Omega$ on $X'(\CC)$.
Using $\Omega$ on $X'(\CC)$, as in Lemma~\ref{lem:upper:estimate:h:0},
we can define $L^2$-norms
of $H^0(X, mL+N)$ and $H^0(X', \pi^*(mL+N))$ as follows: for
$t \in H^0(X, mL+N)$ and
$t' \in H^0(X', \pi^*(mL+N))$,
\[
\Vert t \Vert_{L^2,\Omega}^{m\overline{L}+\overline{N}} =  \left(
\int_{X'(\CC)} \pi^*(\vert t \vert^2_{m\overline{L}+\overline{N}}) \Omega\right)^{1/2}
\]
and
\[
\Vert t' \Vert_{L^2,\Omega}^{\pi^*(m\overline{L}+\overline{N})} = \left(
\int_{X'(\CC)}  \vert t' \vert^2_{\pi^*(m\overline{L}+\overline{N})} \Omega\right)^{1/2},
\]
where $\vert\cdot\vert_{m\overline{L}+\overline{N}}$ and
$\vert\cdot\vert_{\pi^*(m\overline{L}+\overline{N})}$ are  
the hermitian norms of $m\overline{L} + \overline{N}$ and
$\pi^*(m\overline{L}+\overline{N})$ respectively.
Note that $\pi^*(\vert\cdot\vert_{m\overline{L}+\overline{N}}) = \vert\cdot\vert_{\pi^*(m\overline{L}+\overline{N})}$.
Let $\Vert\cdot\Vert_{L^2,\quot}^{\pi^*(m\overline{L}+\overline{N})}$ be
the quotient norm of $\Gamma(X'/X, mL+N)$
induced by $H^0(X', \pi^*(mL+N)) \to \Gamma(X'/X, mL+N)$ and
the $L^2$-norm
$\Vert\cdot\Vert^{\pi^*(m\overline{L}+\overline{N})}_{L^2, \Omega}$
of $H^0(X', \pi^*(mL+N))$. 
Then we have a normed exact sequence 
\addtocounter{Claim}{1}
\begin{multline}
\label{eqn:prop:invariance:vol:birat:morphism:1}
0 \to \left( H^0(X, mL+N), \Vert\cdot\Vert_{L^2}^{m\overline{L}+\overline{N}} \right)
\to \left( H^0(X', \pi^*(mL+N)), \Vert\cdot\Vert_{L^2}^{\pi^*(m\overline{L}+\overline{N})} \right) \\
\to \left( \Gamma(X'/X, mL+N), \Vert\cdot\Vert_{L^2,\quot}^{\pi^*(m\overline{L}+\overline{N})} \right) \to 0. 
\end{multline}
Since $\Vert\cdot\Vert_{L^2,\quot}^{\pi^*(m\overline{L}+\overline{N})}
\leq \Vert\cdot\Vert^{\pi^*(m\overline{L}+\overline{N})}_{\sup,\quot}$, it is sufficient to show that
\[
\ah\left( \Gamma(X'/X, mL+N), \Vert\cdot\Vert^{\pi^*(m\overline{L}+\overline{N})}_{L^2,\quot} \right)
\leq o(m^d).
\]
By virtue of \cite[Corollary~(4.8)]{ZhPL},
$H^0(X', \pi^*(mL+ N))$ is generated by sections $t$ with 
\[
\Vert t \Vert^{\pi^*(m\overline{L}+\overline{N})}_{L^2} \leq \Vert t \Vert^{\pi^*(m\overline{L}+\overline{N})}_{\sup} < 1
\]
for $m \gg 1$ because $\pi^*(\overline{L})$ is ample.
Thus so does
$\Gamma(X'/X, mL+N)$ with respect to $\Vert\cdot\Vert^{\pi^*(m\overline{L}+\overline{N})}_{L^2,\quot}$.
Hence, by using (1) and (5) of Proposition~\ref{prop:basic:property:h:0:h:1:x},  
it suffices to show that
\[
\achi\left( \Gamma(X'/X, mL+N), \Vert\cdot\Vert^{\pi^*(m\overline{L}+\overline{N})}_{L^2,\quot} \right)
\leq o(m^d)
\]
because $\rank \Gamma(X'/X, mL+N) = o(m^{d-1})$.
By using the normed exact sequence \eqref{eqn:prop:invariance:vol:birat:morphism:1} and
\cite[Theorem~2.1, (1)]{Yuan}, we have
\begin{multline*}
\achi\left( \Gamma(X'/X, mL+N), \Vert\cdot\Vert^{\pi^*(m\overline{L}+\overline{N})}_{L^2,\quot} \right) =
\achi\left( H^0(X', \pi^*(mL+N)), \Vert\cdot\Vert_{L^2}^{\pi^*(m\overline{L}+\overline{N})} \right) \\
- \achi\left( H^0(X, mL+N), \Vert\cdot\Vert_{L^2}^{m\overline{L}+\overline{N}} \right) + o(m^d).
\end{multline*}
On the other hand, using \cite[Theorem~(1.4)]{ZhPL} and
Gromov's inequality on $X'(\CC)$, we can see that
\[
\begin{cases}
{\displaystyle \achi\left( H^0(X', \pi^*(mL+N)), \Vert\cdot\Vert_{L^2}^{\pi^*(m\overline{L}+\overline{N})} \right)
= \frac{\adeg(\acherncl_1(\pi^*(\overline{L}))^{\cdot d})}{d!}m^d + o(m^d), } \\
\\
{\displaystyle \achi\left( H^0(X, mL+N), \Vert\cdot\Vert_{L^2}^{m\overline{L}+\overline{N}} \right) =
\frac{\adeg(\acherncl_1(\overline{L})^{\cdot d})}{d!}m^d + o(m^d)}
\end{cases}
\]
as in the proof of Lemma~\ref{lem:upper:estimate:h:0}.
Moreover, by  the projection formula,
\[
\adeg(\acherncl_1(\pi^*(\overline{L}))^{\cdot d})= 
\adeg(\acherncl_1(\overline{L})^{\cdot d}).
\]
Thus the claim follows.

\begin{Claim}
\label{claim:prop:invariance:vol:birat:morphism:2}
If $\pi$ is finite, then
\[
\ah\left( \Gamma(X'/X, mL+N), \Vert\cdot\Vert^{\pi^*(m\overline{L}+\overline{N})}_{\sup,\quot} \right)
\leq o(m^d).
\]
\end{Claim}

Let $\overline{A}$ be an ample $C^{\infty}$-hermitian invertible sheaf on $X$.
Replacing $\overline{A}$ by a higher multiple of $\overline{A}$ if
necessarily, we may assume that there is a non-zero section $s$ of $H^0(X, A - L)$
such that $\Vert s \Vert_{\sup} \leq 1$ and $s$ dose not vanish at any
associated point of $\pi_*(\OO_{X'})/\OO_X$.
Then we have the following commutative diagram:
\[
\begin{CD}
H^0(X', \pi^*(mL+N)) @>{\pi^*(s)}>> H^0(X',  \pi^*(mA+N)) \\
@VVV @VVV \\
\Gamma(X'/X, mL+N) @>{s}>> \Gamma(X'/X, mA+N).
\end{CD}
\]
By our choice of $s$, the horizontal arrows are injective.
Let $\Vert\cdot\Vert^{\pi^*(m\overline{A}+\overline{N})}_{\sup, \pi^*(s), \sub}$ be
the subnorm of $H^0(X', \pi^*(mL+N))$ induced by
\[
H^0(X', \pi^*(mL+N)) \overset{\pi^*(s)}{\longrightarrow} H^0(X', \pi^*(mA+N))
\]
and
$\Vert\cdot\Vert_{\sup}^{\pi^*(m\overline{A}+\overline{N})}$.
Moreover, let $\Vert\cdot\Vert^{\pi^*(m\overline{A}+\overline{N})}_{\sup, \pi^*(s), \sub,\quot}$
be the quotient norm of \\
$\Gamma(X'/X, mL+N)$
induced by
\[
H^0(X', \pi^*(mL+N)) \to \Gamma(X'/X, mL+N),
\]
and let $\Vert\cdot\Vert^{\pi^*(m\overline{A}+\overline{N})}_{\sup,\quot}$
be the quotient norm of $\Gamma(X'/X, mA+N)$  induced by
\[
H^0(X', \pi^*(mA+N)) \to \Gamma(X'/X, mA+N).
\]
Then, by (2) of Lemma~\ref{lem:sub:quot:vs:quot:sub},
\[
\Vert\cdot\Vert^{\pi^*(m\overline{A}+\overline{N})}_{\sup, \pi^*(s), \sub,\quot} \geq
\Vert\cdot\Vert^{\pi^*(m\overline{A}+\overline{N})}_{\sup,\quot} 
\]
on $\Gamma(X'/X, mL+N)$.
Therefore, by the previous claim,
\begin{multline*}
\ah\left( \Gamma(X'/X, mL+N), \Vert\cdot\Vert^{\pi^*(m\overline{A}+\overline{N})}_{\sup, \pi^*(s), \sub,\quot}\right) \\
\leq
\ah\left( \Gamma(X'/X, mA+N), \Vert\cdot\Vert^{\pi^*(m\overline{A}+\overline{N})}_{\sup,\quot} \right) \leq
o(m^d).
\end{multline*}
On the other hand, since 
\[
\Vert\cdot\Vert^{\pi^*(m\overline{L}+\overline{N})}_{\sup,\quot} \geq 
\Vert\cdot\Vert^{\pi^*(m\overline{A}+\overline{N})}_{\sup, \pi^*(s), \sub,\quot},
\]
we have
\begin{multline*}
\ah\left( \Gamma(X'/X, mL+N), \Vert\cdot\Vert^{\pi^*(m\overline{L}+\overline{N})}_{\sup, \quot}\right) \\
\leq
\ah\left( \Gamma(X'/X, mL+N), \Vert\cdot\Vert^{\pi^*(m\overline{A}+\overline{N})}_{\sup, \pi^*(s), \sub,\quot}\right).
\end{multline*}
Thus the claim follows.

\begin{Claim}
\label{claim:prop:invariance:vol:birat:morphism:3}
If $\pi$ is finite, then
\begin{multline*}
\limsup_{m \to \infty} 
\frac{\ah\left(H^0(X, mL + N), \Vert\cdot\Vert_{\sup}^{m\overline{L} + \overline{N}}\right)}{m^d} \\
=
\limsup_{m \to \infty} 
\frac{\ah\left(H^0(X', \pi^*(mL + N)), \Vert\cdot\Vert_{\sup}^{\pi^*(m\overline{L} + \overline{N})}\right)}{m^d}.
\end{multline*}
\end{Claim}

By using (4) of Proposition~\ref{prop:basic:property:h:0:h:1:x} and 
Claim~\ref{claim:prop:invariance:vol:birat:morphism:2},
the normed exact sequence
\begin{multline*}
0 \to \left( H^0(X, mL+N), \Vert\cdot\Vert_{\sup}^{m\overline{L}+\overline{N}} \right)
\to \left( H^0(X', \pi^*(mL+N)), \Vert\cdot\Vert_{\sup}^{\pi^*(m\overline{L}+\overline{N})} \right) \\
\to \left( \Gamma(X'/X, mL+N), \Vert\cdot\Vert_{\sup,\quot}^{\pi^*(m\overline{L}+\overline{N})} \right) \to 0
\end{multline*}
gives rise to
\begin{multline*}
\ah  \left( H^0(X, mL+N), \Vert\cdot\Vert_{\sup}^{m\overline{L}+\overline{N}} \right) \leq
\ah \left( H^0(X', \pi^*(mL+N)), \Vert\cdot\Vert_{\sup}^{\pi^*(m\overline{L}+\overline{N})} \right) \\
\leq \ah  \left( H^0(X, mL+N), \Vert\cdot\Vert_{\sup}^{m\overline{L}+\overline{N}} \right) + o(m^d).
\end{multline*}
This shows the claim.

\medskip
Let us consider a general case.
We set $X'' = \Spec(\pi_*(\OO_{X'}))$.
Then $\pi : X' \to X$ can be factorized
$\pi_1 : X' \to X''$ and $\pi_2 : X'' \to X$ such that
$\pi = \pi_2 \circ \pi_1$,
$(\pi_1)_*(\OO_{X'}) = \OO_{X''}$ and $\pi_2$ is finite.
Thus, by Claim~\ref{claim:prop:invariance:vol:birat:morphism:3},
\begin{multline*}
\limsup_{m \to \infty} 
\frac{\ah\left(H^0(X, mL + N), \Vert\cdot\Vert_{\sup}^{m\overline{L} + \overline{N}}\right)}{m^d} \\
=
\limsup_{m \to \infty} 
\frac{\ah\left(H^0(X'', \pi_2^*(mL + N)), \Vert\cdot\Vert_{\sup}^{\pi_2^*(m\overline{L} + \overline{N})}\right)}{m^d}.
\end{multline*}
On the other hand, since $(\pi_1)_*(\OO_{X'}) = \OO_{X''}$,
\[
H^0(X', \pi^*(mL+N)) = H^0(X'', \pi_2^*(mL+N))
\]
for all $m \geq 1$.
Thus 
\begin{multline*}
\limsup_{m \to \infty} 
\frac{\ah\left(H^0(X'', \pi_2^*(mL + N)), \Vert\cdot\Vert_{\sup}^{\pi_2^*(m\overline{L} + \overline{N})}\right)}{m^d} \\
=
\limsup_{m \to \infty} 
\frac{\ah\left(H^0(X', \pi^*(mL + N)), \Vert\cdot\Vert_{\sup}^{\pi^*(m\overline{L} + \overline{N})}\right)}{m^d}.
\end{multline*}
Hence the theorem follows.
\QED

Next let us consider the following theorem.

\begin{Theorem}
\label{thm:limsup:mL:N}
Let  $\overline{L}$ and $\overline{N}$ be $C^{\infty}$-hermitian invertible sheaves on $X$.
Then
\[
\limsup_{m \to \infty} 
\frac{\ah\left(H^0(X, mL + N), \Vert\cdot\Vert_{\sup}^{m\overline{L} + \overline{N}}\right)}{m^d}
=
\frac{\avol(\overline{L})}{d!}.
\]
\end{Theorem}

\Proof
By Theorem~\ref{thm:invariance:vol:birat:morphism},
we may assume that $X$ is generically smooth.
By using Lemma~\ref{lem:ample:generic:smooth},
there are ample $C^{\infty}$-hermitian invertible sheaves $\overline{A}$ and $\overline{B}$
such that $-\overline{B} \leq \overline{N} \leq \overline{A}$ and that
$\overline{A}$ and $\overline{B}$ satisfy the assumptions (i), (ii) and (iii) of
Theorem~\ref{thm:h:0:estimate:big}.
The inequalities $-\overline{B} \leq \overline{N} \leq \overline{A}$ gives rise to
\begin{multline*}
 \limsup_{m \to \infty} 
\frac{\ah\left(H^0(X, mL -B ), \Vert\cdot\Vert_{\sup}^{m\overline{L} - \overline{B}}\right)}{m^d} \\
\leq
 \limsup_{m \to \infty} 
\frac{\ah\left(H^0(X, mL + N), \Vert\cdot\Vert_{\sup}^{m\overline{L} + \overline{N}}\right)}{m^d} \\
\leq
 \limsup_{m \to \infty} 
\frac{\ah\left(H^0(X, mL + A), \Vert\cdot\Vert_{\sup}^{m\overline{L} + \overline{A}}\right)}{m^d}.
\end{multline*}
Applying Theorem~\ref{thm:h:0:estimate:big} to the case where $b=1$ and $c=0$,
we have
\[
\ah\left(H^0(X, mL + A), \Vert\cdot\Vert_{\sup}^{m\overline{L} + \overline{A}}\right) \leq
\ah\left(H^0(X, mL), \Vert\cdot\Vert_{\sup}^{m\overline{L}}\right) + o(m^d)
\]
for $m \gg 1$, which yields
\[
 \limsup_{m \to \infty} 
\frac{\ah\left(H^0(X, mL + A), \Vert\cdot\Vert_{\sup}^{m\overline{L} + \overline{A}}\right)}{m^d}
\leq \frac{\avol(\overline{L})}{d!}.
\]
Further, 
applying Theorem~\ref{thm:h:0:estimate:big} to the case where $b=c=1$,
\[
\ah\left(H^0(X, mL ), \Vert\cdot\Vert_{\sup}^{m\overline{L}}\right) \leq
\ah\left(H^0(X, mL -B ), \Vert\cdot\Vert_{\sup}^{m\overline{L} - \overline{B}}\right) + o(m^d)
\]
for $m \gg 1$, which implies
\[
\frac{\avol(\overline{L})}{d!} \leq
 \limsup_{m \to \infty} 
\frac{\ah\left(H^0(X, mL -B ), \Vert\cdot\Vert_{\sup}^{m\overline{L} - \overline{B}}\right)}{m^d}.
\]
Thus we get the theorem.
\QED

The following lemma is need to see the characterization of bigness and
the homogeneity of the arithmetic volume function.

\begin{Lemma}
\label{lem:big:property}
Let $\overline{L}$ and $\overline{N}$ be $C^{\infty}$-hermitian invertible sheaves on $X$.
We assume that $\overline{L}$ is big.
Then, for a fixed positive integer $p$,
\[
\limsup_{n\to\infty} \frac{\ah(H^0(X, pnL+N), \Vert\cdot\Vert_{\sup}^{pn\overline{L}+\overline{N}})}{(pn)^d} =
\limsup_{m\to\infty}\frac{\ah(H^0(X, mL+N), \Vert\cdot\Vert_{\sup}^{m\overline{L}+\overline{N}})}{m^d}
\]
and
\[
\liminf_{n\to\infty}
\frac{\ah(H^0(X, pnL+N), \Vert\cdot\Vert_{\sup}^{pn\overline{L}+\overline{N}})}{(pn)^d} =
\liminf_{m\to\infty}
\frac{\ah(H^0(X, mL+N), \Vert\cdot\Vert_{\sup}^{m\overline{L}+\overline{N}})}{m^d}
\]
\end{Lemma}

\Proof
First we claim the following:

\begin{Claim}
\label{claim:lem:big:property:1}
There is a positive integer $m_0$ such that
$\ah(H^0(X, mL), \Vert\cdot\Vert_{\sup}^{m \overline{L}}) \not= 0$
for all $m \geq m_0$.
\end{Claim}

Let $\overline{A}$ be an ample $C^{\infty}$-hermitian invertible sheaf on $X$
such that 
\[
\ah(H^0(X, A), \Vert\cdot\Vert_{\sup}^{\overline{A}}) \not= 0\quad\text{and}\quad
\ah(H^0(X, L + A), \Vert\cdot\Vert_{\sup}^{\overline{L} + \overline{A}}) \not= 0.
\]
Since $\overline{L}$ is big,
we can find a positive integer $a$
with $\ah(H^0(X, aL - A), \Vert\cdot\Vert_{\sup}^{a \overline{L} - \overline{A}}) \not= 0$
(cf. \cite[Proposition~2.2]{MoArHt}).
Note that
\[
aL = (aL - A) + A\quad\text{and}\quad
(a+1)L = (aL - A) + (L + A).
\]
Thus
\[
\ah(H^0(X, aL), \Vert\cdot\Vert_{\sup}^{a \overline{L}}) \not= 0\quad\text{and}\quad
\ah(H^0(X, (a+1)L), \Vert\cdot\Vert_{\sup}^{(a+1) \overline{L}}) \not= 0.
\]
Let $m$ be an integer with $m \geq a^2 + a$.
We set $m = aq + r$, where $0 \leq r < a$.
Then $q \geq a$. Thus we can find $b > 0$ with $q = b + r$.
Therefore
$m \overline{L} = b(a \overline{L}) + r ((a+1)\overline{L})$, which means that
\[
\ah(H^0(X, mL), \Vert\cdot\Vert_{\sup}^{m \overline{L}}) \not= 0.
\]

\medskip
Next we claim the following:

\begin{Claim}
\label{claim:lem:big:property:2}
There is a positive integer $n_0$ such that
\begin{multline*}
\qquad
\ah(H^0(X, pnL+N), \Vert\cdot\Vert_{\sup}^{pn\overline{L}+\overline{N}}) \\
\leq  \ah(H^0(X, (p(n+n_0) + i)L+N), \Vert\cdot\Vert_{\sup}^{(p(n+n_0) + i)\overline{L}+\overline{N}}) \\
\leq \ah(H^0(X, p(n+2n_0 +1)L+N), \Vert\cdot\Vert_{\sup}^{p(n+2n_0 +1)\overline{L}+\overline{N}})
\end{multline*}
for all $n \geq 1$ and all $i = 0, \ldots, p$.
\end{Claim}

We choose $n_0$ with $pn_0 \geq m_0$.
For each $i=0, \ldots, p$, there is a non-zero section
$s_i$ of $H^0(X, (pn_0 + i)L)$  with $\Vert s_i \Vert_{\sup} \leq 1$.
Therefore we have injective homomorphisms
\[
H^0(X, pnL+N) \overset{s_i}{\longrightarrow} H^0(X, (p(n+n_0)+i)L+N)
\overset{s_{p-i}}{\longrightarrow} H^0(p(n + 2n_0 + 1)L+N).
\]
Thus our the claim follows.

\medskip
Let us go back to the proof of the lemma.
By the above claim,
\begin{multline*}
\limsup_{n \to \infty}
\frac{\ah(H^0(X, pnL+N), \Vert\cdot\Vert_{\sup}^{pn\overline{L}+\overline{N}})}{(pn)^d} \\
\leq  \limsup_{n \to \infty}
\frac{\ah(H^0(X, (p(n+n_0) + i)L), \Vert\cdot\Vert_{\sup}^{(p(n+n_0) + i)\overline{L}+\overline{N}})}{(pn)^d} \\
\leq \limsup_{n \to \infty}
\frac{\ah(H^0(X, p(n+2n_0 +1)L), \Vert\cdot\Vert_{\sup}^{p(n+2n_0 +1)\overline{L}+\overline{N}})}{(pn)^d}.
\end{multline*}
Note that
\[
\lim_{n \to \infty} \frac{(pn)^d}{(p(n+n_0)+i)^d} =
\lim_{n \to \infty} \frac{(pn)^d}{(p(n+2n_0 + 1))^d} = 1.
\]
This shows that
\begin{multline*}
\limsup_{n \to \infty}
\frac{\ah(H^0(X, pnL+N), \Vert\cdot\Vert_{\sup}^{pn\overline{L}+\overline{N}})}{(pn)^d} \\
=
\limsup_{n \to \infty}
\frac{\ah(H^0(X, (pn+ i)L+N), \Vert\cdot\Vert_{\sup}^{(pn+i)\overline{L}+\overline{N}})}{(pn+ i)^d} \\
\end{multline*}
for all $i=0, \ldots, p-1$.
Hence
\[
\limsup_{n\to\infty} \frac{\ah(H^0(X, pnL+N), \Vert\cdot\Vert_{\sup}^{pn\overline{L}+\overline{N}})}{(pn)^d} =
\limsup_{m\to\infty}\frac{\ah(H^0(X, mL+N), \Vert\cdot\Vert_{\sup}^{m\overline{L}+\overline{N}})}{m^d}.
\]
In the same way, we can see
\[
\liminf_{n\to\infty}
\frac{\ah(H^0(X, pnL+N), \Vert\cdot\Vert_{\sup}^{pn\overline{L}+\overline{N}})}{(pn)^d} =
\liminf_{m\to\infty}
\frac{\ah(H^0(X, mL), \Vert\cdot\Vert_{\sup}^{m\overline{L}+\overline{N}})}{m^d}.
\]
\QED

The following theorem is a characterization of a big
$C^{\infty}$-hermitian invertible sheaf.
The similar property is observed in \cite{Yuan}.

\begin{Theorem}
\label{thm:equiv:big}
For a $C^{\infty}$-hermitian invertible sheaf $\overline{L}$ on $X$,
the following are equivalent:
\begin{enumerate}
\renewcommand{\labelenumi}{(\arabic{enumi})}
\item
$\avol(\overline{L}) > 0$.

\item $\overline{L}$ is big.

\item
${\displaystyle \liminf_{m\to\infty} \frac{\ah(H^0(X, mL), \Vert\cdot\Vert_{\sup})}{m^d} > 0}$.

\item
${\displaystyle 
\liminf_{m\to\infty} \frac{ \log \# \{ s \in H^0(X, mL) \mid \Vert s \Vert_{\sup} < 1 \}}{m^d} > 0}$.
\end{enumerate}
\end{Theorem}

\Proof
Obviously (3) $\Longrightarrow$ (1) and
(4) $\Longrightarrow$ (1), so that it is sufficient to show
that (1) $\Longrightarrow$ (2), (2) $\Longrightarrow$ (3) and
(2) $\Longrightarrow$ (4).

\medskip
(1) $\Longrightarrow$ (2):
We assume that $\avol(\overline{L}) > 0$.
By (3) of Proposition~\ref{prop:basic:properties:volume},
there is a positive integer $m$ and
a non-zero section $s$ of $H^0(X, mL)$ with
$\Vert s \Vert_{\sup}^{m \overline{L}} < 1$.
Let $\overline{A}$ be an ample $C^{\infty}$-hermitian invertible sheaf on $X$.
By Theorem~\ref{thm:limsup:mL:N},
\[
\limsup_{m \to \infty} 
\frac{\ah\left(H^0(X, mL -A), \Vert\cdot\Vert_{\sup}^{m\overline{L} - \overline{A}}\right)}{m^d}
=
\frac{\avol(\overline{L})}{d!} > 0,
\]
which implies that there is a positive integer $n$
with $\ah\left(H^0(X, nL -A), \Vert\cdot\Vert_{\sup}^{n\overline{L} - \overline{A}}\right) \not= 0$.
Hence $n\overline{L} \geq \overline{A}$.
In particular, $L_{\QQ}$ is big on $X_{\QQ}$.

\medskip
(2) $\Longrightarrow$ (3):
Let $\overline{A}$ be an ample $C^{\infty}$-hermitian invertible sheaf on $X$.
Since $\overline{L}$ is big, there is a positive integer $p$ with
$p\overline{L} \geq \overline{A}$.
Therefore,
\[
\liminf_{n\to\infty} \frac{\ah(H^0(pnL), \Vert\cdot\Vert_{\sup}^{pn\overline{L}})}{(pn)^d}
\geq \frac{1}{p^d} \liminf_{n\to\infty} \frac{\ah(H^0(nA), \Vert\cdot\Vert_{\sup}^{n\overline{A}})}{n^d} > 0.
\]
Hence, by Lemma~\ref{lem:big:property},
\[
\liminf_{m\to\infty} \frac{\ah(H^0(mL), \Vert\cdot\Vert_{\sup}^{m\overline{L}})}{m^d} > 0.
\]

\medskip
(2) $\Longrightarrow$ (4):
We choose a sufficiently small positive number $\lambda$ such that
$\overline{L}^{-\lambda}$ is big.
Since (2) $\Longrightarrow$ (3), we have
\[
\liminf_{m\to\infty}
 \frac{ \log \# \{ s \in H^0(X, mL) \mid \exp(m\lambda) \Vert s \Vert_{\sup} \leq  1 \}}{m^d} > 0,
\]
which yields (4).
\QED

\begin{Remark}
In the paper \cite{Yuan},
Yuan uses the condition (4) of the above theorem
as a definition of a big $C^{\infty}$-hermitian invertible sheaf.
By the above theorem, Yuan's definition is equivalent to our bigness.
\end{Remark}

\begin{Proposition}
\label{prop:avol:hom}
$\avol$ is homogeneous of degree $d$, that is,
$\avol(p\overline{L}) = p^d \avol(\overline{L})$ for every non-negative
integer $p$.
\end{Proposition}

\Proof
Since
\[
\limsup_{n \to \infty}
\frac{\ah(H^0(X, npL), \Vert\cdot\Vert_{\sup}^{np \overline{L}})}{(np)^d} \leq
\limsup_{m \to \infty}
\frac{\ah(H^0(X, mL), \Vert\cdot\Vert_{\sup}^{m \overline{L}})}{m^d},
\]
we have $\avol(p \overline{L}) \leq p^d \avol(\overline{L})$.
Thus, if $\avol(\overline{L}) = 0$, then the assertion is obvious.
Therefore we may assume that $\avol(\overline{L}) > 0$,
namely, by Theorem~\ref{thm:equiv:big}, $\overline{L}$ is big.
Hence, by Lemma~\ref{lem:big:property},
\[
\limsup_{n\to\infty} \frac{\ah(H^0(X, npL), \Vert\cdot\Vert_{\sup}^{pn\overline{L}})}{(np)^d} =
\limsup_{m\to\infty}\frac{\ah(H^0(X, mL), \Vert\cdot\Vert_{\sup}^{m\overline{L}})}{m^d},
\]
which means that $\avol(p\overline{L}) = p^d \avol(\overline{L})$.
\QED


\section{Continuity of the volume function}
Let $X$ be a $d$-dimensional projective arithmetic variety
and $\aPic(X)$ the group of isomorphism classes of
$C^{\infty}$-hermitian invertible sheaves on $X$.
An element of $\aPic(X) \otimes \QQ$ is called
a {\em $C^{\infty}$-hermitian $\QQ$-invertible sheaf} on $X$.
For $\overline{L} \in \aPic(X)$, 
the image of $\overline{L}$ via
the canonical homomorphism
$\aPic(X) \to \aPic(X) \otimes \QQ$ is denoted by $[\overline{L}]$.
Note that $[\overline{L}] = [(\OO_X, \vert\cdot\vert_{can})]$ if and only if
$\overline{L}$ is a torsion in $\aPic(X)$, that is,
there is a positive integer $n$ with
$n \overline{L} = (\OO_X, \vert\cdot\vert_{can})$.
We say 
a $C^{\infty}$-hermitian $\QQ$-invertible sheaf
$\overline{L}$ is represented by $\overline{M} \in \aPic(X)$
if $[\overline{M}] = \overline{L}$.
Moreover a $C^{\infty}$-hermitian $\QQ$-invertible sheaf $\overline{L}$
on $X$ is said to be ample if
there is a positive integer $n$ such that
$n\overline{L}$ is represented by
an ample $C^{\infty}$-hermitian invertible sheaf on $X$.
Similarly we say $\overline{L}$ is nef (resp. big) if
$n \overline{L}$ is represented by
a nef (resp. big) $C^{\infty}$-hermitian invertible sheaf for some positive integer $n$.
Let us begin with the following lemma.

\begin{Lemma}
\label{lem:avol:Pic:Q}
$\avol : \aPic(X) \to \RR$ extends to a homogeneous map
\[
\avol: \aPic(X) \otimes \QQ \to \RR
\]
of degree $d$, that is, $\avol(a\overline{L}) = a^d \avol(\overline{L})$ for
every non-negative rational number $a$.
\end{Lemma}

\Proof
Let $\overline{L}$ be a $C^{\infty}$-hermitian $\QQ$-invertible sheaf on $X$.
Let $n$ be a positive integer such that
$n \overline{L}$ is represented by a $C^{\infty}$-hermitian invertible sheaf $\overline{M}$.
Then we would like to define
$\avol(\overline{L})$ to be
$\avol(\overline{M})/n^d$.
Indeed this is well-defined.
Let $n'$ be another positive integer
such that
$n' \overline{L}$ is represented by a $C^{\infty}$-hermitian invertible sheaf $\overline{M}'$.
Then, since $[n'\overline{M}] = [n\overline{M}']$,
there is a positive integer $m$ with
$m n'\overline{M} = m n\overline{M}'$. On the other hand,
\[
\avol(mn'\overline{M}) = (mn')^d \avol(\overline{M})\quad\text{and}\quad
\avol(mn\overline{M}') = (mn)^d \avol(\overline{M}')
\]
Thus $\avol(\overline{M})/n^d = \avol(\overline{M}')/{n'}^d$.

Next let us see that $\avol(a\overline{L}) = a^d \avol(\overline{L})$
for every non-negative rational number $a$.
Let $n$ and $m$ be positive integers such that
$m a \in \ZZ$ and $n\overline{L}$ is represented by $\overline{M} \in \aPic(X)$.
Then, since $(mn)a\overline{L}$ is represented by $(ma) \overline{M}$,
\[
\avol(a\overline{L}) = \avol((ma) \overline{M})/(mn)^d = a^d \avol(\overline{M})/n^d = 
a^d \avol(\overline{L}).
\]
\QED

In Conventions and terminology~\ref{CT:Arakelov:order:hermitian},
we define the order $\leq$ on the group $\aPic(X)$.
We would like to extend it to 
$\aPic(X) \otimes \QQ$.
For $\overline{L}, \overline{M} \in \aPic(X) \otimes \QQ$,
if there is a positive integer $n$
such that
$n \overline{L}$ and $n \overline{M}$ are represented by
a $C^{\infty}$-hermitian invertible sheaf $\overline{L}'$ and
$\overline{M}'$ respectively with $\overline{L}' \leq \overline{M}'$, then
we denote this by $\overline{L} \leq_{\QQ} \overline{M}$.

\begin{Lemma}
\label{lem:L<M:vL<vM}
For $\overline{L}, \overline{L}', \overline{M}, \overline{M}' \in \aPic(X) \otimes \QQ$,
we have the following:
\begin{enumerate}
\renewcommand{\labelenumi}{(\arabic{enumi})}
\item
$\overline{L} \leq_{\QQ} \overline{M}$ if and only if
$-\overline{M} \leq_{\QQ} -\overline{L}$.

\item
If $\overline{L} \leq_{\QQ} \overline{M}$ and
$\overline{L}' \leq_{\QQ} \overline{M}'$, then
$\overline{L} + \overline{L}' \leq_{\QQ} \overline{M} + \overline{M}'$.

\item
If $\overline{L} \leq_{\QQ} \overline{M}$ and $a$ is a non-negative rational number,
then $a \overline{L} \leq_{\QQ} a \overline{M}$.

\item
If $\overline{L} \leq_{\QQ} \overline{M}$, then $\avol(\overline{L}) \leq \avol(\overline{M)}$.
\end{enumerate}
\end{Lemma}

\Proof
(1), (2) and (3) are consequence of the properties in
Conventions and terminology~\ref{CT:Arakelov:order:hermitian}.
Let us consider (4).
Let $n$ be a positive integer
such that $n\overline{L}$ and $n\overline{M}$ are
represented by 
$C^{\infty}$-hermitian invertible sheaves $\overline{L}'$ and
$\overline{M}'$ with $\overline{L}' \leq \overline{M}'$.
Then $\avol(\overline{L}') \leq \avol(\overline{M}')$
by (1) of Proposition~\ref{prop:basic:properties:volume}.
Hence we have (4).
\QED

\begin{Remark}
\label{rem:order:hermitian:vector:space}
For reader's convenience, let us give a sketch of the proof of
the properties (1) and (2) in Conventions and terminology~\ref{CT:Arakelov:order:hermitian}.
Let $(V, \sigma)$ and $(W, \tau)$ be normed $\CC$-vector spaces of dimension one.
We denote $(V, \sigma) \leq (W, \tau)$ if
there is an isomorphism $\phi : V \to W$ over $\CC$ such that
$\tau(\phi(x)) \leq \sigma(x)$ for all $x \in V$.
Then, in order to see the properties (1) and (2), it is sufficient to show the following:
\begin{enumerate}
\renewcommand{\labelenumi}{(\alph{enumi})}
\item
$(V, \sigma) \leq (W, \tau)$ if and only if
$(W^{\vee}, \tau^{\vee}) \leq (V^{\vee}, \sigma^{\vee})$.

\item
If $(V, \sigma) \leq (W, \tau)$ and $(V', \sigma') \leq (W', \tau')$, then
\[
(V \otimes V', \sigma \otimes \sigma') \leq (W \otimes W', \tau \otimes \tau').
\]
\end{enumerate}

\medskip
(a) 
Let $\phi : V \to W$ be an isomorphism over $\CC$,
$v$ a basis of $V$ and $w = \phi(v)$.
Let $v^{\vee}$ and $w^{\vee}$ be the dual bases of $v$ and $w$ respectively.
Since $\sigma(v/\sigma(v)) = 1$,
\[
\sigma^{\vee}(v^{\vee}) = \max \{ \vert v^{\vee}(x) \vert \mid \sigma(x) = 1 \}
= 1/\sigma(v).
\]
In the same way, $\tau^{\vee}(w^{\vee}) = 1/\tau(w)$.
Note that $\phi^{\vee}(w^{\vee}) = v^{\vee}$.
Thus (a) follows.

\medskip
(b) Let $\phi : V \to W$ and $\phi' : V' \to W'$ be isomorphisms over $\CC$
such that $\tau(\phi(x)) \leq \sigma(x)$ and $\tau'(\phi'(x')) \leq \sigma'(x')$
for all $x \in V$ and $x' \in V'$.
Then
\[
(\tau \otimes \tau')((\phi \otimes \phi')(x \otimes x')) =
\tau(\phi(x)) \tau'(\phi'(x')) \leq \sigma(x) \sigma'(x') = (\sigma \otimes \sigma')(x \otimes x').
\]
Therefore $(V \otimes V', \sigma \otimes \sigma') \leq (W \otimes W', \tau \otimes \tau')$.
\end{Remark}

The following theorem is the main result of this paper.

\begin{Theorem}[Continuity of volume]
\label{thm:cont:arithmetic:volume}
Let $\overline{L}$ and $\overline{A}$ be
$C^{\infty}$-hermitian $\QQ$-invertible sheaves on $X$.
Then
\[
\lim_{\substack{\epsilon \in \QQ \\ \epsilon \to 0}}
\avol(\overline{L} +\epsilon \overline{A}) = \avol(\overline{L}).
\]
More generally,
for $C^{\infty}$-hermitian $\QQ$-invertible sheaves $\overline{A}_1, \ldots ,\overline{A}_n$
on $X$,
\[
\lim_{\substack{\epsilon_1, \ldots, \epsilon_n \in \QQ \\ \epsilon _1\to 0, \ldots , \epsilon_n \to 0}}
\avol(\overline{L} +\epsilon_1 \overline{A}_1 + \cdots + \epsilon_n \overline{A}_n ) = \avol(\overline{L}).
\]
\end{Theorem}

\Proof
First let us consider the case $n=1$.
Let $\mu : X' \to X$ be a generic resolution of singularities of $X$.
Then, by Theorem~\ref{thm:invariance:vol:birat:morphism},
$\avol(\overline{L} +\epsilon \overline{A}) = \avol(\mu^*(\overline{L}) +\epsilon \mu^*(\overline{A}))$
and
$\avol(\overline{L}) = \avol(\mu^*(\overline{L}))$.
Thus we may assume that $X$ is generically smooth.

\begin{Claim}
\label{claim:thm:cont:arithmetic:volume:1}
We may further assume that $\overline{A}$ is ample.
\end{Claim}

Let $\overline{B}$ be an ample $C^{\infty}$-hermitian $\QQ$-invertible sheaf on $X$
such that $\overline{A} + \overline{B}$ is ample.
Then, for $\epsilon \geq 0$,
\[
\overline{L} - \epsilon(\overline{A} + \overline{B}) \leq_{\QQ} \overline{L} - \epsilon \overline{A}
\leq_{\QQ} \overline{L} + \epsilon \overline{B}
\quad\text{and}\quad
\overline{L}  - \epsilon \overline{B}
\leq_{\QQ}
\overline{L} + \epsilon \overline{A}
\leq_{\QQ}
\overline{L} + \epsilon(\overline{A} + \overline{B}).
\]
Thus, by (4) of Lemma~\ref{lem:L<M:vL<vM},
\[
\begin{cases}
\avol(\overline{L} - \epsilon(\overline{A} + \overline{B})) \leq 
\avol(\overline{L} - \epsilon \overline{A})
\leq \avol(\overline{L} + \epsilon \overline{B}), \\
\avol(\overline{L}  - \epsilon \overline{B})
\leq
\avol(\overline{L} + \epsilon \overline{A})
\leq
\avol(\overline{L} + \epsilon(\overline{A} + \overline{B})).
\end{cases}
\]
Hence the claim follows.

\medskip
From now on, we assume that $\overline{A}$ is ample.
It is obvious that 
\[
\lim_{\substack{\epsilon \in \QQ \\ \epsilon \to 0}}
\avol(\overline{L} +\epsilon \overline{A}) = \avol(\overline{L})
\quad\Longleftrightarrow\quad
\lim_{\substack{\epsilon \in \QQ \\ \epsilon \to 0}}
\avol(\overline{L} +\epsilon a \overline{A}) = \avol(\overline{L})
\]
for any positive rational number $a$.
Moreover,
\[
\avol(n\overline{L} + \epsilon \overline{A}) = n^d \vol(\overline{L} + (\epsilon/n) \overline{A})
\quad\text{and}\quad
\avol(n\overline{L}) = n^d \vol(\overline{L} ).
\]
Therefore, 
we may assume that 
$\overline{L}$ is  $C^{\infty}$-hermitian invertible sheaf.
Further, by Lemma~\ref{lem:ample:generic:smooth},
we may assume that
$\overline{A}$ is a $C^{\infty}$-hermitian invertible sheaf and that
$\overline{A}$ satisfies the assumptions (i), (ii) and (iii) of
Theorem~\ref{thm:h:0:estimate:big}.

Since
\[
\avol(\overline{L} - \epsilon' \overline{A}) \leq
\avol(\overline{L} - \epsilon \overline{A}) \leq 
\avol(\overline{L}) \leq
\avol(\overline{L} + \epsilon \overline{A}) \leq
\avol(\overline{L} - \epsilon' \overline{A})
\]
for $0 \leq \epsilon \leq \epsilon'$, it is sufficient to show that
\[
\avol(\overline{L}) = \lim_{p \to \infty} \avol(\overline{L} + (1/p)\overline{A})
= \lim_{p \to \infty} \avol(\overline{L} - (1/p)\overline{A}).
\]
By Theorem~\ref{thm:h:0:estimate:big}
(or Proposition~\ref{prop:0:estimate:big:curve} for $d=1$),
there are positive constants $a'_0$, $C'$ and $D'$
depending only on $X$, $\overline{L}$ and $\overline{A}$ such that
\begin{multline*}
\ah\left(H^0(a L + (b-c)A),
\Vert\cdot\Vert^{a \overline{L} + (b-c) \overline{A}}_{\sup}\right) \leq
\ah\left(H^0(a L - cA), 
\Vert\cdot\Vert^{a \overline{L} - c \overline{A}}_{\sup}\right) \\
+ C' b a^{d-1} + D' a^{d-1} \log(a)
\end{multline*}
for all integers $a, b, c$ with $a \geq b \geq c \geq 0$ and $a \geq a'_0$.

First we set $a = pm$, $b = m$ and $c = 0$ for a fixed positive integer $p$. Then
\begin{multline*}
\ah\left(H^0(pm L + m A),
\Vert\cdot\Vert^{pm \overline{L} + m \overline{A}}_{\sup}\right) \leq
\ah\left(H^0(pm  L), 
\Vert\cdot\Vert^{pm \overline{L}}_{\sup}\right) \\
+ C' p^{d-1} m^{d} + D' p^{d-1} m^{d-1} \log(pm)
\end{multline*}
for $m \gg 1$.
This implies that
\[
\avol(p \overline{L} + \overline{A}) \leq
\avol(p \overline{L}) + C' p^{d-1}
\]
for all $p \geq 1$, which means that
\[
\avol(\overline{L}) \leq \avol( \overline{L} + (1/p)\overline{A}) \leq
\avol(\overline{L}) + C'(1/p).
\]
Hence
\[
\lim_{p \to \infty}
\avol( \overline{L} + (1/p)\overline{A}) = \avol(\overline{L}).
\]

Next we set $a = pm$ and $b = c = m$.
Then
\begin{multline*}
\ah\left(H^0(pm L),
\Vert\cdot\Vert^{pm \overline{L}}_{\sup}\right) \leq
\ah\left(H^0(pm L - mA), 
\Vert\cdot\Vert^{pm \overline{L} - m \overline{A}}_{\sup}\right) \\
+ C' p^{d-1} m^{d} + D' p^{d-1} m^{d-1}\log(pm)
\end{multline*}
for $m \gg 1$. This implies that
\[
\avol(\overline{L} - (1/p)\overline{A}) \leq \avol(\overline{L}) \leq
\avol(\overline{L} - (1/p)\overline{A}) + C' (1/p).
\]
Thus
\[
\lim_{p \to \infty}
\avol( \overline{L} - (1/p)\overline{A}) = \avol(\overline{L}).
\]

\medskip
Let us consider a general case.
We can find $C^{\infty}$-hermitian $\QQ$-invertible sheaves $\overline{A}'_i$
and $\overline{A}''_i$ such that $0 \leq_{\QQ} \overline{A}'_i$,
$0 \leq_{\QQ} \overline{A}''_i$ and
$\overline{A}_i = \overline{A}'_i - \overline{A}''_i$ for each $i$.
Then
\[
\overline{L}+  \epsilon_1 \overline{A}_1 + \cdots + \epsilon_n \overline{A}_n =
\overline{L}+  \epsilon_1 \overline{A}'_1 + \cdots + \epsilon_n \overline{A}'_n +
(- \epsilon_1) \overline{A}''_1 + \cdots + (-\epsilon_n) \overline{A}''_n.
\]
Thus we may
assume that $0 \leq_{\QQ} \overline{A}_1, \ldots, 
0 \leq_{\QQ} \overline{A}_n$.
Find an ample $C^{\infty}$-hermitian $\QQ$-invertible sheaf $\overline{B}$ such that
$\overline{A}_i \leq_{\QQ} \overline{B}$ for all $i=1, \ldots, n$. Then
\[
-\vert \epsilon_i \vert B \leq_{\QQ} -\vert \epsilon_i \vert \overline{A}_i \leq_{\QQ}
 \epsilon_i \overline{A}_i \leq_{\QQ}
\vert \epsilon_i \vert \overline{A}_i \leq_{\QQ} \vert \epsilon_i \vert \overline{B}
\]
for each $i$, which implies
\[
\overline{L} -(\vert \epsilon_1 \vert + \cdots + \vert \epsilon_n \vert ) \overline{B} \leq_{\QQ}
\overline{L} + \epsilon_1 \overline{A}_1 + \cdots + \epsilon_n \overline{A}_n \leq_{\QQ}
\overline{L} + (\vert \epsilon_1 \vert + \cdots + \vert \epsilon_n \vert ) \overline{B}.
\]
Therefore
\begin{multline*}
\avol(\overline{L} -(\vert \epsilon_1 \vert + \cdots + \vert \epsilon_n \vert ) \overline{B}) \\
\leq
\avol(\overline{L}+  \epsilon_1 \overline{A}_1 + \cdots + \epsilon_n \overline{A}_n) \\
\leq
\avol(\overline{L} + (\vert \epsilon_1 \vert + \cdots + \vert \epsilon_n \vert ) \overline{B}).
\end{multline*}
Thus the general assertion follows from the case $n=1$.
\QED

As a corollary, we can show the following arithmetic Hilbert-Samuel theorem for
a nef $C^{\infty}$-hermitian invertible sheaf.

\begin{Corollary}[Arithmetic Hilbert-Samuel formula]
\label{cor:Hilbert:Samuel:nef}
Let $\overline{L}$ and $\overline{N}$ 
be $C^{\infty}$-hermitian invertible sheaves on $X$.
If $\overline{L}$ is nef, then
\[
\ah\left(H^0(X, mL+N), \Vert\cdot\Vert_{\sup}^{m \overline{L}+\overline{N}}\right) = 
\frac{\adeg (\acherncl_1(\overline{L})^{\cdot d})}{d!}{m^d} + o(m^d)\quad(m \gg 1).
\]
In particular, $\avol(\overline{L}) = \adeg (\acherncl_1(\overline{L})^{\cdot d})$,
and $\overline{L}$ is big if and only if
$ \adeg (\acherncl_1(\overline{L})^{\cdot d}) > 0$.
\end{Corollary}

\Proof
First let us see the following claim:

\begin{Claim}
\label{claim:cor:Hilbert:Samuel:nef:1}
$\avol(\overline{L}) = \adeg (\acherncl_1(\overline{L})^{\cdot d})$.
\end{Claim}

Let $\overline{A}$ be an ample $C^{\infty}$-hermitian
invertible sheaf on $X$.
Then $\overline{L} + \epsilon \overline{A}$ is ample for all $\epsilon > 0$.
Thus
\[
\avol(\overline{L} + \epsilon \overline{A}) =
\adeg \left( \left( \acherncl_1(\overline{L}) + \epsilon \acherncl_1(\overline{A}) \right)^{\cdot d}
\right).
\]
Therefore our claim follows from the continuity of volumes.

\medskip
Let us go back to the proof of the corollary.
It is sufficient to show
\[
\frac{\adeg (\acherncl_1(\overline{L})^{\cdot d})}{d! }=\lim_{m \to \infty}
\frac{\ah(H^0(X, mL+N), \Vert\cdot\Vert_{\sup}^{m \overline{L}+\overline{N}})}{m^d}.
\]
If $\overline{L}$ is not big, then,
by  Claim~\ref{claim:cor:Hilbert:Samuel:nef:1},
\[
\avol(\overline{L}) = \adeg (\acherncl_1(\overline{L})^{\cdot d}) = 0.
\]
Thus our assertion is obvious by Theorem~\ref{thm:limsup:mL:N}, 
so that
we may assume that $\overline{L}$ is big.
Then there is a positive integer $k$ with $k\overline{L} \geq \overline{A}$.
We set $\overline{E} = k \overline{L} - \overline{A}$.
Since
\[
p \overline{L} - \overline{E} = (p - k) \overline{L} + \overline{A},
\]
$p \overline{L} - \overline{E}$ is ample if $p \geq k$. 
On the other hand, since $p \overline{L} \geq p \overline{L} - \overline{E}$,
we have
\[
\ah(H^0(X, npL+N), \Vert\cdot\Vert_{\sup}^{np \overline{L}+\overline{N}}) \geq
\ah(H^0(X, n(pL-E)+N), \Vert\cdot\Vert_{\sup}^{n(p \overline{L} - \overline{E})+\overline{N}})
\]
for $n \geq 1$, which implies that
\begin{multline*}
\liminf_{n \to \infty}
\frac{\ah(H^0(X, npL+N), \Vert\cdot\Vert_{\sup}^{np \overline{L}+\overline{N}})}{(np)^d} \\
\geq
\liminf_{n \to \infty}
\frac{\ah(H^0(X, n(pL-E)+N), \Vert\cdot\Vert_{\sup}^{n(p \overline{L} - \overline{E})+\overline{N}})}{(np)^d}.
\end{multline*}
Therefore, for a fixed $p$ with $p \geq k$, by using  Lemma~\ref{lem:upper:estimate:h:0} and
Lemma~\ref{lem:big:property},
\[
\liminf_{m \to \infty}
\frac{\ah(H^0(X, mL+N), \Vert\cdot\Vert_{\sup}^{m \overline{L}+\overline{N}})}{m^d} \geq
\frac{\adeg (\acherncl_1(p \overline{L} - \overline{E})^{\cdot d})}{p^d d!}.
\]
Thus, taking $p \to \infty$,
\[
\liminf_{m \to \infty}
\frac{\ah(H^0(X, mL+N), \Vert\cdot\Vert_{\sup}^{m \overline{L}+\overline{N}})}{m^d} \geq
\frac{\adeg (\acherncl_1(\overline{L})^{\cdot d})}{d!}.
\]
On the other hand, by Theorem~\ref{thm:limsup:mL:N} and Claim~\ref{claim:cor:Hilbert:Samuel:nef:1},
\[
\limsup_{m \to \infty}
\frac{\ah(H^0(X, mL+N), \Vert\cdot\Vert_{\sup}^{m \overline{L}+\overline{N}})}{m^d}
= \frac{\adeg (\acherncl_1(\overline{L})^{\cdot d})}{d!},
\]
which proves the corollary.
\QED

Finally let us consider the volume of the difference of nef
$C^{\infty}$-hermitian $\QQ$-invertible sheaves, which is essentially
the main result of Yuan's paper \cite{Yuan}.

\begin{Theorem}
\label{thm:vol:diff:nef}
Let $\overline{L}$ and $\overline{M}$ be nef $C^{\infty}$-hermitian $\QQ$-invertible sheaves on $X$.
Then
\[
\avol(\overline{L} - \overline{M}) \geq \adeg(\acherncl_1(\overline{L})^{\cdot d}) -
d \cdot \adeg(\acherncl_1(\overline{L})^{\cdot (d-1)} \cdot \acherncl_1(\overline{M})).
\]
\end{Theorem}

\Proof
First we assume that $\overline{L}$ and $\overline{M}$ are
ample $C^{\infty}$-hermitian invertible sheaves on $X$.
Then, by \cite{Yuan},
\begin{multline*}
\avol(\overline{L} - \overline{M}) \geq
\limsup_{m\to\infty}
\frac{\chi\left( H^0(m(L - M)), \Vert\cdot\Vert^{m(\overline{L}-\overline{M})}_{\sup} \right)}{m^d/d!} \\
\geq \adeg(\acherncl_1(\overline{L})^{\cdot d}) -
d \cdot \adeg(\acherncl_1(\overline{L})^{\cdot (d-1)} \cdot \acherncl_1(\overline{M})).
\end{multline*}
Thus, using the homogeneity of $\avol$,
the inequality holds for ample $C^{\infty}$-hermitian $\QQ$-invertible sheaves on $X$.
Let $\overline{A}$ be an ample $C^{\infty}$-hermitian invertible sheaf on $X$.
Then, for a small positive number $\epsilon$,
$\overline{L} + \epsilon \overline{A}$ and
$\overline{M} + \epsilon \overline{A}$ are ample. Thus,
\begin{multline*}
\avol(\overline{L} - \overline{M}) = \avol((\overline{L} + \epsilon \overline{A}) - (\overline{M} + \epsilon \overline{A})) \\
\geq 
\adeg((\acherncl_1(\overline{L})+ \epsilon \acherncl_1(\overline{A}))^{\cdot d}) \\
-d \cdot \adeg((\acherncl_1(\overline{L})+ \epsilon \acherncl_1(\overline{A}))^{\cdot (d-1)} \cdot (\acherncl_1(\overline{M})+ \epsilon \acherncl_1(\overline{A}))).
\end{multline*}
Therefore the theorem follows.
\QED

\begin{Remark}[Arithmetic analogue of Fujita's approximation theorem]
It is very natural to ask the following arithmetic analogue of Fujita's approximation theorem:
Let $\overline{L}$ be a big $C^{\infty}$-hermitian $\QQ$-invertible sheaf on $X$.
For any positive number $\epsilon$, do there exist a birational morphism
$\mu : X' \to X$ and an ample $C^{\infty}$-hermitian $\QQ$-invertible sheaf $\overline{A}$ on $X'$
such that $\overline{A} \leq_{\QQ} \mu^*(\overline{L})$ and
$\avol(\overline{L}) \leq \avol(\overline{A}) + \epsilon$ ?
\end{Remark}

\section{Generalized Hodge index theorem}
In this section, we consider a generalized Hodge index theorem
as an application of the continuity of the volume function.
First let us introduce a technical definition.

Let $X$ be a projective arithmetic variety of dimension $d$.
Let $L$ be an invertible sheaf on $X$ such that
$L$ is nef on the generic fiber $X_{\QQ}$ of $X \to \Spec(\ZZ)$.
We say $L$ has {\em moderate growth of positive even cohomologies} if
there are a generic resolution of singularities $\mu : Y \to X$ and
an ample invertible sheaf $A$ on $Y$ such that,
for any positive integer $n$, there is a positive integer $m_0$ such that
\[
\log \# (H^{2i}(Y, m(n\mu^*(L) + A))) = o(m^d)
\]
for all $m \geq m_0$ and for all $i > 0$.
Here we consider examples of invertible sheaves with moderate growth of
positive even cohomologies.

\begin{Example}
\label{example:moderate:growth:coho}
(1) We assume that $d = 2$. Then $L$ has obviously moderate growth of positive even cohomologies.

\medskip
(2) If $L$ is nef on each geometric fiber of $X \to \Spec(\ZZ)$, then
$L$ has moderate growth of positive even cohomologies.
Indeed, let $\mu : Y \to X$ be a generic resolution of singularities and
$A$ an ample invertible sheaf on $Y$. Then, for all $n \geq 1$,
$n \mu^*(L) + A$ is ample. Thus $H^{i}(Y, m(n\mu^*(L) + A)) = 0$
for $m \gg 1$ and $i > 0$.

\medskip
(3) We assume that $d = 2$ and $X$ is generically smooth.
Let $E$ be a rank $r$ locally free sheaf on $X$.
Let $\pi : P = \Proj(\bigoplus_{m\geq 0} \Sym^m(E)) \to X$ be the projective bundle of $E$ and
$\OO_P(1)$ the tautological invertible sheaf of $P$.
We set $L = r \cdot \OO_P(1) - \pi^*(\det E)$.
If $E$ is semistable on the generic fiber $X_{\QQ}$, then
it is well-known that $L$ is nef on the generic fiber $P_{\QQ}$.
Moreover $L$ has moderate growth of positive even cohomologies. 
This fact can be checked as follows:
Let $B$ be an ample invertible sheaf on $X$ such that
$A = \OO_P(1) + \pi^*(B)$ is ample. Then
\begin{multline*}
H^i(P, m(nL + A)) \\
= H^i(P,  \OO_P(mnr + m) + \pi^*(mB - mn \det(E))) \\
= H^i(X, \Sym^{mnr + m}(E) \otimes (mB - mn \det(E)))
\end{multline*}
because $R^j\pi_* \OO_P(l) = 0$ for $l \geq 0$ and $j > 0$.
In particular, 
\[
H^i(P, m(nL + A)) = 0
\]
for $i \geq 2$.
\end{Example}

The main result of this section is the following generalized Hodge index theorem.

\begin{Theorem}[Generalized Hodge index theorem]
\label{thm:ineq:avol:deg}
Let $X$ be a $d$-dimensional projective arithmetic variety and $\overline{L}$ a
$C^{\infty}$-hermitian invertible sheaf on $X$.
We assume the following:
\begin{enumerate}
\renewcommand{\labelenumi}{(\arabic{enumi})}
\item
$L_{\QQ}$ is nef on $X_{\QQ}$.

\item
$c_1(\overline{L})$ is semipositive on $X(\CC)$.

\item
$L$ has moderate growth of positive even cohomologies.
\end{enumerate}
Then we have
an inequality $\avol(\overline{L}) \geq \adeg(\acherncl_1(\overline{L})^{\cdot d})$.
\end{Theorem}

\Proof
First we assume that $X$ is generically smooth. Moreover,
instead of the properties (1), (2) and (3) as above, we assume the following
\rom{(a)}, \rom{(b)} and \rom{(c)}:
\begin{enumerate}
\renewcommand{\labelenumi}{(\alph{enumi})}
\item
$L_{\QQ}$ is ample on $X_{\QQ}$.

\item
$c_1(\overline{L})$ is positive on $X(\CC)$.

\item
There is a positive number $m_0$ such that
\[
\log \#(H^{2i}(X, mL)) = o(m^d)
\]
for $m \geq m_0$.
\end{enumerate}
Then let us see the following:

\begin{Claim}
\label{claim:thm:ineq:avol:deg:1}
$\avol(\overline{L}) \geq \adeg(\acherncl_1(\overline{L})^{\cdot d})$.
\end{Claim}

By virtue of the arithmetic Riemann-Roch theorem \cite{GSArRR} and the asymptotic estimate of
analytic torsions \cite{BV}, we obtain
\begin{multline*}
\achi(H^0(X, mL), \Vert\cdot\Vert_{L^2}^{\overline{L}}) \\
+\sum_{i \geq 1}
\log \#(H^{2i}(X, mL)) - \sum_{i \geq 1}  \log \#(H^{2i-1}(X, mL)) \\
= \frac{\adeg(\acherncl_1(\overline{L})^{\cdot d})}{d!} m^d + o(m^d)
\end{multline*}
for $m \gg 1$. Thus, using the assumption (c) and (1) of 
Proposition~\ref{prop:basic:property:h:0:h:1:x},
\[
\ah(H^0(X, mL), \Vert\cdot\Vert_{L^2}^{\overline{L}}) \geq
\frac{\adeg(\acherncl_1(\overline{L})^{\cdot d})}{d!} m^d + o(m^d)
\]
for $m \gg 1$. By Gromov's inequality and (3) of
Proposition~\ref{prop:basic:property:h:0:h:1:x}, the above inequality implies
\[
\ah(H^0(X, mL), \Vert\cdot\Vert_{\sup}^{\overline{L}}) \geq
\frac{\adeg(\acherncl_1(\overline{L})^{\cdot d})}{d!} m^d + o(m^d)
\]
for $m \gg 1$.  Hence the claim follows.

\medskip
Let us go back to a general case.
Since $L$ has moderate growth of positive even cohomologies,
there is a generic resolution of singularities $\mu : Y \to X$ and
an ample invertible sheaf $A$ on $Y$ such that,
for any positive integer $n$, there is a positive integer $m_0$ such that
$\log \# (H^{2i}(Y, m(n\mu^*(L) + A))) = o(m^d)$ for all $m \geq m_0$ and
for all $i > 0$.
Let us give a $C^{\infty}$-hermitian metric $\vert \cdot\vert_A$ to $A$
such that $\overline{A} = (A, \vert\cdot\vert_A)$ is ample as
a $C^{\infty}$-hermitian invertible sheaf.
Then, by Claim~\ref{claim:thm:ineq:avol:deg:1},
\[
\avol(n \mu^*(\overline{L}) + \overline{A}) \geq \adeg(\acherncl_1(n\mu^*(\overline{L}) + \overline{A})^{\cdot d}),
\]
which implies 
\[
\avol(\mu^*(\overline{L}) + (1/n)\overline{A}) \geq \adeg((\acherncl_1(\mu^*(\overline{L})) + (1/n)\acherncl_1(\overline{A}))^{\cdot d})
\]
by Proposition~\ref{prop:avol:hom}.
Hence, using the continuity of the volume function,
\[
\avol(\mu^*(\overline{L})) \geq \adeg(\acherncl_1(\mu^*(\overline{L}))^{\cdot d}).
\]
This gives rise to our assertion by Theorem~\ref{thm:invariance:vol:birat:morphism}
and the projection formula.
\QED

According to (1), (2) and (3) of Example~\ref{example:moderate:growth:coho},
we have the following corollaries.

\begin{Corollary}
\label{cor:HIT:dim:2}
Let $X$ be a projective arithmetic surface and $\overline{L}$ a $C^{\infty}$-hermitian invertible
sheaf on $X$ such that
$L$ is nef on the generic fiber of $X \to \Spec(\ZZ)$ and
$c_1(\overline{L})$ is semipositive on $X(\CC)$.
Then 
\[
\avol(\overline{L}) \geq \adeg(\acherncl_1(\overline{L})^{\cdot 2}).
\]
\end{Corollary}

\begin{Corollary}
\label{cor:HIT:nef}
Let $X$ be a projective arithmetic variety of dimension $d$ and $\overline{L}$ a $C^{\infty}$-hermitian invertible
sheaf on $X$ such that
$L$ is nef on every geometric fiber of $X \to \Spec(\ZZ)$ and
$c_1(\overline{L})$ is semipositive on $X(\CC)$.
Then 
\[
\avol(\overline{L}) \geq \adeg(\acherncl_1(\overline{L})^{\cdot d}).
\]
In particular, if $\adeg(\acherncl_1(\overline{L})^{\cdot d}) > 0$, then
$\overline{L}$ is big.
\end{Corollary}

\begin{Corollary}
\label{cor:Bogo:dim:2}
Let $X$ be a projective and generically smooth
arithmetic surface and $\overline{E}$ a $C^{\infty}$-hermitian locally
free sheaf on $X$.
If  the metric of $\overline{E}$ is
Einstein-Hermitian, then
\[
\adeg\left( \acherncl_2(\overline{E}) - \frac{r-1}{2r} \acherncl_1(\overline{E})^2 \right) \geq 0,
\]
where $r = \rank E$.
\end{Corollary}

\Proof
Let $\pi : P = \Proj(\bigoplus_{m\geq 0} \Sym^m(E)) \to X$ be the projective bundle of $E$ and
$\OO_P(1)$ the tautological invertible sheaf of $P$.
Using the surjective homomorphism $\pi^*(E) \to \OO_P(1)$ and the hermitian metric of $E$,
we give the quotient metric $\vert\cdot\vert_P$ to $\OO_P(1)$.
We set $\overline{\OO}_P(1) = (\OO_P(1), \vert\cdot\vert_P)$ and
$\overline{L} = r \cdot \overline{\OO}_P(1) - \pi^*(\det \overline{E})$.
Note that $L_{\QQ}$ is nef and not big and that
$c_1(\overline{L})$ is semipositive (cf. \cite[Lemma~8.7.1]{MoBogo1}). 
Moreover $L$ has moderate growth of positive even cohomologies. 
Thus, by Theorem~\ref{thm:ineq:avol:deg},
\[
\adeg(\acherncl_1(\overline{L})^{\cdot r+1}) \leq 0
\]
because $\overline{L}$ is not big.
Note that
\[
\adeg(\acherncl_1(\overline{L})^{\cdot r+1}) = r^{r+1} \cdot
\adeg\left( \frac{r-1}{2r} \acherncl_1(\overline{E})^2 - \acherncl_2(\overline{E}) \right)
\]
(cf. \cite[Section~8]{MoBogo1}).
Thus
\[
\adeg\left( \acherncl_2(\overline{E}) - \frac{r-1}{2r} \acherncl_1(\overline{E})^2 \right) \geq 0.
\]
\QED

\begin{Remark}
(1) In Corollary~\ref{cor:HIT:dim:2}, if $\deg(L_{\QQ}) = 0$,
then $\adeg(\acherncl_1(\overline{L})^{\cdot 2}) \leq 0$.
This is nothing more than the Hodge index theorem due to Faltings and
Hriljac (cf. \cite{Fal} and \cite{Hri}). In this sense, we call Theorem~\ref{thm:ineq:avol:deg}
the generalized Hodge index theorem.

\medskip
(2) The second assertion of Corollary~\ref{cor:HIT:nef} 
is a generalization of \cite[Corollary~(1.9)]{ZhPL}.

\medskip
(3) Corollary~\ref{cor:Bogo:dim:2} is valid even if $E_{\QQ}$ is semistable on $X_{\QQ}$.
The case where the metric is Einstein-Hermitian is however essential and crucial for
a general case.
For details, see \cite{MoBogo1}.
\end{Remark}

\end{document}